\documentclass[11pt]{article}

\usepackage{amssymb,latexsym}
\usepackage{amsmath,amscd}
\usepackage{theorem}
\usepackage[all]{xy}

\title{\bf Dirac structures, momentum maps and quasi-Poisson manifolds}
\author{Henrique Bursztyn\thanks{{\tt henrique@math.toronto.edu}}\\[0.1cm]
        Department of Mathematics\\
    University of Toronto\\
        Toronto, Ontario M5S 3G3, Canada
        \\[0.2cm]
        Marius Crainic \thanks{{\tt crainic@math.uu.nl}} 
        \\[0.1cm]
         Department of Mathematics\\
        Utrecht University, P.O. Box 80.010, 3508 TA\\
        Utrecht, The Netherlands}
\date{}

\newcommand{\st} [1]     {\scriptscriptstyle{{#1}}}
\newcommand{\phiG}        {\phi^{{G}}}
\newcommand{\tri}        {\chi}
\newcommand{\rmap}{\longrightarrow}

\newcommand{\id}         {{\mathrm {Id}}}

\newcommand{\Ker}        {{\mathrm {Ker}}}

\newcommand{\Ad}         {{\mathrm {Ad}}}
\newcommand{\gra}        {{\mathrm {graph}}}
\newcommand{\SP} [1]     {{\left\langle {{#1}} \right\rangle}}
\newcommand{\Bi} [1]     {{\left( {{#1}} \right)_{\mathfrak{g}}}}

\newcommand{\pr}         {{\mathrm{pr}}}

\newcommand{\Jac}       {\mathrm{Jac}}

\newcommand{\piQ}      {\pi_{{Q}}}
\newcommand{\piP}      {\pi_{{P}}}
\newcommand{\LN}       {L_{{N}}}
\newcommand{\LM}       {L_{{M}}}

\newcommand{\grd}         {\mathcal{G}}
\newcommand{\sour}        {\mathsf{s}}
\newcommand{\tar}         {{\mathsf{t}}}

\newcommand{\Cour}[1]      {[\![#1]\!]}

\newcommand{\Lie}        {\mathcal L}

\newtheorem{lemma} {Lemma} [section]
\newtheorem{proposition} [lemma] {Proposition}
\newtheorem{claim} [lemma] {Claim}
\newtheorem{theorem} [lemma] {Theorem}
\newtheorem{corollary} [lemma] {Corollary}
\newtheorem{definition}[lemma] {Definition}

\theorembodyfont{\rm} 

\newtheorem{example}[lemma] {Example}

\newtheorem{remark}[lemma]{Remark}

\newenvironment{proof}{{\sc Proof:}}{{\hspace*{\fill} $\square$\\}}

\numberwithin{equation}{section}

\begin{document}

\maketitle

\begin{abstract}
We extend the correspondence between Poisson maps and actions of
symplectic groupoids, which generalizes the one between momentum
maps and hamiltonian actions, to the realm of Dirac geometry. As
an example, we show how hamiltonian quasi-Poisson manifolds fit
into this framework by constructing  an ``inversion'' procedure
relating quasi-Poisson bivectors to twisted Dirac structures.
\end{abstract}

\begin{center}
{\bf Dedicated to Alan Weinstein for his 60th birthday}
\end{center}


\tableofcontents

\section{Introduction}

This paper builds on three ideas pursued by Alan
Weinstein in some of his many fundamental contributions to Poisson geometry: First,
 Lie algebroids play a prominent role in the study of Poisson manifolds \cite{SilWein99,We98};
second, Poisson maps can be regarded as generalized momentum maps
for actions of symplectic groupoids \cite{MiWe,We02}; third,
Poisson structures on manifolds are particular examples of more
general objects, called Dirac structures
\cite{Cou90,CouWe,SeWe01}. The main objective of this paper is to
combine these three ideas in order to extend the notion of
``momentum map'' to the realm of Dirac geometry.
As an application, we obtain an alternative approach to
hamiltonian quasi-Poisson manifolds \cite{AKM} which answers many
of the questions posed in \cite{SeWe01,We02}, shedding light on
the relationship between various notions of generalized Poisson
structures, hamiltonian actions and reduced spaces.

Let $\mathfrak{g}$ be a Lie algebra, and consider its dual
$\mathfrak{g}^*$, equipped with its Lie-Poisson structure. The
central ingredients in the formulation of classical hamiltonian
$\mathfrak{g}$-actions are a Poisson manifold $(Q,\piQ)$ and a
Poisson map $J:Q \to \mathfrak{g}^*$, which we use to define an
action of $\mathfrak{g}$ on $Q$ by hamiltonian vector fields:
\begin{equation}\label{eq:hamiltaction}
\mathfrak{g} \longrightarrow \mathcal{X}(Q),\;\;\; v \mapsto
X_{J_v}:= i_{dJ_v}(\piQ),
\end{equation}
where $J_v \in C^\infty(Q)$ is given by $J_v(x)= \SP{J(x),v}$. For
the global picture, we assume that $J$ is a \textit{complete}
Poisson map \cite[Sec.~6.2]{SilWein99}, in which case the
infinitesimal action \eqref{eq:hamiltaction} can be integrated to
an action of the connected, simply connected Lie group $G$ with
Lie algebra $\mathfrak{g}$, in such a way that $J$ becomes
$G$-equivariant with respect to the coadjoint action of $G$ on
$\mathfrak{g}^*$. The map $J$ is called a \textbf{momentum map}
for the $G$-action on $Q$, and we refer to the $G$-action as
\textbf{hamiltonian}.
A key observation, described in \cite{MiWe,We02}, is that this
construction of a hamiltonian action out of a Poisson map  holds
in much more generality: one may replace $\mathfrak{g}^*$ by any
Poisson manifold, as long as Lie groups are replaced by symplectic
groupoids \cite{We87}. In this sense, any Poisson map can be seen
as a ``Poisson-manifold valued moment map''.

In this paper, we show that the correspondence between Poisson maps
and hamiltonian actions by symplectic groupoids can be further
extended to the context  of Dirac geometry: in this setting,
Poisson maps must be replaced by special types of Dirac
maps, called  Dirac realizations (see Def.~\ref{def:diracreal}); for
the associated global actions, twisted presymplectic groupoids \cite{BCWZ}
(alternatively called quasi-symplectic groupoids \cite{Xu}) play
the role of symplectic groupoids. Our main results show that various important
notions of generalized  hamiltonian actions, such as the ``quasi'' objects of
\cite{AKM,AMM}, fit nicely into the Dirac geometry framework.

We organize our results as follows:

In Section \ref{sec:liealg}, we discuss important connections between Lie algebroids and
bivector fields. Our main result is that, just as ordinary Poisson structures give rise to
Lie algebroid structures on their cotangent bundles, a quasi-Poisson manifold
\cite[Def.~2.1]{AKM} $(M,\pi)$  defines a Lie algebroid structure on
$T^*M\oplus \mathfrak{g}$, where $\mathfrak{g}$ is the Lie algebra of the Lie group acting
on $M$. 
The leaves of this Lie algebroid coincide with the leaves of the
``quasi-hamiltonian foliation'' of \cite[Sec.~9]{AKM} in the hamiltonian case, though,
in our framework, we make no assumption about the existence of group-valued moment maps.

In Section \ref{sec:infinit}, we study  hamiltonian actions in the
context of Dirac geometry at the infinitesimal level. We observe
that, just as Poisson maps, Dirac realizations are always
associated with Lie algebroid actions. (This is, in fact, the
guiding principle in our definition of Dirac realizations.) After
discussing how classical notions of infinitesimal hamiltonian
actions fit into this framework, we prove the main result of the
section: Dirac realizations of Cartan-Dirac structures on Lie
groups \cite{BCWZ,SeWe01} are equivalent to 
quasi-Poisson $\mathfrak{g}$-manifolds carrying group-valued moment maps. 
This equivalence involves an ``inversion'' procedure relating twisted Dirac structures and
quasi-Poisson bivectors, revealing that these two objects are in a
certain sense ``mirror'' to one another. The main ingredients in
this discussion
 are the Lie algebroids of Section \ref{sec:liealg}
and the  bundle maps  which appear in \cite{BCWZ} as infinitesimal
versions of multiplicative 2-forms. This result explains, in
particular, the relationship between Cartan-Dirac and
quasi-Poisson structures on Lie groups; on the other hand, it
recovers the correspondence proven in \cite[Thm.~10.3]{AKM}
between ``non-degenerate'' hamiltonian quasi-Poisson manifolds
(i.e., those for which the Lie algebroids of Section
\ref{sec:liealg} are \textit{transitive}) and
 quasi-hamiltonian spaces \cite{AMM}.

In Section \ref{sec:global}, we study moment maps in Dirac geometry from
a global point of view. We show that \textit{complete} Dirac realizations ``integrate''
to presymplectic groupoid actions, which are natural extensions of those studied in
\cite{Xu}. As our main example,
we show that the  ``integration'' of Dirac realizations of Cartan-Dirac structures
on Lie groups results in hamiltonian quasi-Poisson $G$-manifolds.
Finally, we show that the natural reduction procedure in the setting of Dirac geometry
 encompasses various classical reduction theorems
\cite{Lu,MW,MiWe} as well as their ``quasi'' counterparts \cite{AKM,AMM,Xu}.

We remark, following an observation of E. Meinrenken, that the results concerning
quasi-Poisson manifolds in this paper only require the Lie algebras to be quadratic,
in contrast with some
of the constructions in \cite{AKM}, in which the positivity of the bilinear forms plays a key role.
(In particular, our results hold for quasi-Poisson $G$-manifolds when $G$
is a noncompact semisimple Lie group.) Most of our constructions can be carried out
in the more general setting of \cite{AK}, but this will be discussed in a separate paper.

A work which gave initial motivation  and is closely related to
the present paper is that of Xu \cite{Xu}, in which a Morita
theory of quasi-symplectic groupoids is developed in order to
compare ``moment map theories''. Our results show that twisted
Dirac structures complement Xu's picture in two ways: on one hand,
by providing the infinitesimal framework for Morita equivalence;
on the other hand, by leading to more general ``modules''
(i.e., hamiltonian spaces).

It is a pleasure to dedicate this paper to Alan Weinstein, whose work and insightful ideas
have been an unlimited source of inspiration to us.

\smallskip

\noindent \textbf{Acknowledgments:} We would like to thank many
people for helpful discussions concerning this work, including A.
Alekseev, Y. Kosmann-Schwarzbach, E. Meinrenken, D. Roytenberg, A. Weinstein,
P. Xu and the referees. Our collaboration was facilitated by invitations to the
conferences ``PQR2003'', in Brussels, and ``AlanFest'' and
``Symplectic Geometry and Moment Maps'', held at the Erwin
Schr\"odinger Institute, where most of the results in this paper
were announced; we thank the organizers of these meetings, in
particular A. Alekseev, S. Gutt, J. Koiller, J. Marsden, and T. Ratiu. 
For financial support, H.B. thanks DAAD (German Academic
Exchange Service) and M.C. thanks  KNAW (Dutch Royal Academy of Arts and Sciences).
H.B. thanks Freiburg University for its hospitality while
part of this work was being done.

\smallskip

\noindent \textbf{Notation:} We use the following conventions for bundle maps:
if $\pi$ is a bivector field on $M$, then ${\pi}^{\sharp}:T^*M\to TM$, $\alpha \mapsto
\pi(\alpha,\cdot)$; if $\omega$ is a 2-form, then ${\omega}^{\sharp}:TM \to T^*M$, $X
\mapsto \omega(X,\cdot)$. 

The space of $k$-multivector fields on $M$ is denoted by $\mathcal{X}^k(M)$.

On a Lie group $G$, with Lie algebra $\mathfrak{g}$, 
$\Bi{\cdot,\cdot}$ will denote a bi-invariant nondegenerate quadratic form;
we write $\phiG$ for the associated Cartan 3-form, and $\tri_G\in \Lambda^3 \mathfrak{g}$
for the dual trivector. The Lie algebra $\mathfrak{g}$ is identified with \textit{right}-invariant
vector fields on $G$.

\section{Lie algebroids, bivector fields and Poisson geometry}
\label{sec:liealg}

\subsection{Lie algebroids}

A \textbf{Lie algebroid} over a manifold $M$ is a vector bundle
$A\to M$ together with a Lie algebra bracket $[\cdot,\cdot]$ on
the space of sections $\Gamma(A)$, and a bundle map $\rho: A \to
TM$, called the \textbf{anchor}, satisfying the Leibniz identity
\begin{equation}\label{eq:leibniz}
[\xi,f \xi']=f[\xi,\xi']+\Lie_{\rho(\xi)}(f)\xi', \;\;\; \mbox{
for } \xi, \xi' \in \Gamma(A), \; \mbox{ and } f \in C^\infty(M).
\end{equation}
Whenever there is no risk of confusion, we will write
$\Lie_{\rho(\xi)}$ simply as $\Lie_{\xi}$.

If $M$ is a point, then a Lie algebroid over $M$ is a Lie algebra
in the usual sense. An important feature of Lie algebroids $A \to
M$ is that the image of the anchor, $\rho(A) \subseteq TM$,
defines a generalized integrable distribution, determining a
singular foliation of $M$. The leaves of this foliation are the
\textbf{orbits} of the Lie algebroid. The following example plays
a key role in the study of hamiltonian actions and moment maps.

\begin{example}(\textit{Transformation Lie algebroids})\label{ex:transfalg}

Consider an infinitesimal action of a  Lie algebra $\mathfrak{g}$
on a manifold $M$, given by a Lie algebra homomorphism
$\bar{\rho}: \mathfrak{g} \to \mathcal{X}(M)$. The
\textbf{transformation Lie algebroid} associated with this action
is the trivial vector bundle $M\times \mathfrak{g}$, with anchor
$(x,v) \mapsto \rho(x,v):= \bar{\rho}(v)(x)$ and Lie bracket on
$\Gamma(M\times \mathfrak{g})=C^\infty(M,\mathfrak{g})$ defined by
\begin{equation}\label{eq:actionbracket}
[u,v](x) := [u(x),v(x)]_{\mathfrak{g}} + (\bar{\rho}(u(x))\cdot v)(x) -
 (\bar{\rho}(v(x))\cdot u)(x).
\end{equation}
We often denote a transformation Lie algebroid by $\mathfrak{g}\ltimes M$.

Note that $[\cdot,\cdot]$ is uniquely determined by the condition that it coincides
with $[\cdot,\cdot]_{\mathfrak{g}}$ on constant functions and the Leibniz identity.
The orbits of $\mathfrak{g}\ltimes M$ are the $\mathfrak{g}$-orbits on $M$.

\end{example}

The remaining of this section is devoted to examples of Lie
algebroids closely related to Poisson manifolds.

\subsection{Bivector fields and Poisson structures}\label{subsec:liealgpoiss}

If $(M,\pi)$ is a Poisson manifold, then $T^*M$ carries a Lie
algebroid structure with anchor
\begin{equation}\label{eq:pisharp}
\pi^{\sharp}:T^*M \to TM,\;\;\;
\beta(\pi^{\sharp}(\alpha))=\pi(\alpha,\beta),
\end{equation}
and bracket
\begin{equation}\label{eq:bracket1}
 [\alpha, \beta]= \Lie_{\pi^{\sharp}(\alpha)}(\beta)-
\Lie_{\pi^{\sharp}(\beta)}(\alpha)- d\pi(\alpha, \beta),
\end{equation}
uniquely characterized by $[df, dg]= d\{f, g\}$ and the Leibniz
identity \eqref{eq:leibniz}.
Here, as usual, $\{f, g\} = \pi(df, dg)$ is the Poisson bracket on
$C^{\infty}(M)$. In this case, the orbits of $T^*M$ are the
symplectic leaves of $M$, i.e., the integral manifolds of the
distribution defined by the hamiltonian vector fields $X_f=
\pi^{\sharp}(df)$.

\begin{example} \textit{(Lie-Poisson structures)} \label{Lie-Poisson}

Let $(\mathfrak{g},[\cdot,\cdot])$ be a Lie algebra, and consider
$\mathfrak{g}^*$ equipped with the associated Lie-Poisson
structure
\begin{equation}\label{eq:liepoiss}
\{f,g\}(\mu):=\SP{\mu,[df(\mu),dg(\mu)]},\;\; \mu \in
\mathfrak{g}^*.
\end{equation}
Under the identification $T^*\mathfrak{g}^*\cong \mathfrak{g}^*\times \mathfrak{g}$,
one can see that the Lie algebroid structure on $T^*\mathfrak{g}^*$ induced by
\eqref{eq:liepoiss} is that of a transformation Lie algebroid
$\mathfrak{g}\ltimes \mathfrak{g}^*$, see Example
\ref{ex:transfalg}, and a  direct computation reveals that the
action in question is the coadjoint action.
\end{example}

If $\pi \in \mathcal{X}^2(M)$ is an \textit{arbitrary} bivector
field, let us consider $\pi^{\sharp}$, $[\cdot, \cdot]$, $\{\cdot,
\cdot\}$ and $X_{f}$ as defined by the previous formulas, and let
$\tri_{\pi} \in \mathcal{X}^3(M)$ be the trivector field defined
by
\begin{equation}
\tri_{\pi}:=[\pi,\pi],
\end{equation}
i.e., $\tri_{\pi}$ satisfies
\[
\frac{1}{2}\tri_{\pi}(df,dg,dh) = \{ f, \{g, h\}\}+ \{ g, \{h, f\}\}+ \{
h, \{f, g\}\} = \{f,\{g,h\}\} + c.p.,
\]
where we use $c.p.$ to denote cyclic permutations.
\begin{lemma}
\label{bivector} For any bivector field $\pi$ on $M$, one has
\begin{align}
&\pi^{\sharp}([\alpha, \beta])= [\pi^{\sharp}(\alpha),
\pi^{\sharp}(\beta)]- \frac{1}{2}i_{\alpha\wedge\beta}(\tri_\pi) ,\label{eq:ident1}\\
&[\alpha,[\beta, \gamma]] + c.p. =
\frac{1}{2}(\Lie_{i_{\alpha\wedge\beta}(\tri_{\pi})}(\gamma) + c.p.)-
d(\tri_{\pi}(\alpha, \beta, \gamma)),\label{eq:ident2}
\end{align}
for $\alpha, \beta, \gamma \in \Omega^1(M)$. As a result, the
following are equivalent:
\begin{enumerate}
\item[(i)] $\pi$ is a Poisson tensor;

\item[(ii)] $\pi^{\sharp}:\Omega^1(M) \to \mathcal{X}(M)$
preserves the brackets;

\item[(iii)] the bracket $[\cdot, \cdot]$ on $\Omega^1(M)$
satisfies the Jacobi identity;

\item[(iv)] $(T^*M, \pi^{\sharp}, [\cdot, \cdot])$ is a Lie
algebroid.
\end{enumerate}
\end{lemma}

\begin{proof} The key remark is that the difference between the left and right hand
sides of each of \eqref{eq:ident1} and \eqref{eq:ident2} is
$C^{\infty}(M)$-multilinear in $\alpha$, $\beta$ and $\gamma$. So
it is enough to prove the identities on exact forms, which is
immediate.
\end{proof}

\begin{example}(\textit{Twisted Poisson manifolds})\label{ex:twispoiss}

Consider a closed 3-form $\phi\in \Omega^3(M)$. A
\textbf{$\phi$-twisted Poisson structure} on $M$
\cite{Strobl,Park} consists of a bivector field $\pi \in
\mathcal{X}^2(M)$ satisfying
\[
\frac{1}{2}[\pi, \pi]= \pi^{\sharp}(\phi).
\]
Here, we abuse notation and write $\pi^{\sharp}$ to denote the map
induced by \eqref{eq:pisharp} on exterior algebras. We know from
Lemma \ref{bivector} that the bracket \eqref{eq:bracket1} induced
by $\pi$ is not preserved by $\pi^{\sharp}$ and does not satisfy
the Jacobi identity. However,
\[
\pi^{\sharp}([\alpha, \beta]+ i_{\pi^{\sharp}(\alpha)\wedge
\pi^{\sharp}(\beta)}(\phi))= [\pi^{\sharp}(\alpha),
\pi^{\sharp}(\beta)] .
\]
Hence, if we define a ``twisted'' version of the bracket
\eqref{eq:bracket1},
\[
[\alpha, \beta]_{\phi}:= [\alpha, \beta]+
i_{\pi^{\sharp}(\alpha)\wedge \pi^{\sharp}(\beta)}(\phi),
\]
then $\pi^{\sharp}$ will preserve this new bracket, and
$[\cdot,\cdot]_{\phi}$ satisfies the Jacobi identity. As a result,
$(T^*M, \pi^{\sharp}, [\cdot, \cdot]_{\phi})$ is a Lie algebroid.
We leave it to the reader to prove a ``twisted'' version of Lemma
\ref{bivector}.
\end{example}

\subsection{The Lie algebroid of a quasi-Poisson manifold}\label{subsec:liealgq}
Let $G$ be a Lie group with Lie algebra $\mathfrak{g}$, equipped
with a bi-invariant nondegenerate quadratic form $\Bi{\cdot, \cdot}$. Let $\phiG$
be the bi-invariant Cartan 3-form on $G$, and let $\tri_G \in \Lambda^3\mathfrak{g}$
be its dual trivector. On Lie algebra elements $u,v, w \in \mathfrak{g}$, we have
$$
\phiG(u,v,w) = \tri_G(u^{\vee},v^{\vee},w^{\vee}) = \frac{1}{2}\Bi{u,[v,w]},
$$
where $u^{\vee},v^{\vee},w^{\vee}$ are dual to $u,v,w$ via $\Bi{\cdot,\cdot}$;
when $\Bi{\cdot,\cdot}$ is a metric and $e_a$ is an orthonormal basis of $\mathfrak{g}$, 
we can write\footnote{More generally, with no positivity assumptions on $\Bi{\cdot,\cdot}$, 
we can write $\tri_{G}= \frac{1}{12} \sum \Bi{e_a, [e_b, e_c]} f_a\wedge f_b\wedge f_c$,
where $f_a$ is a basis of $\mathfrak{g}$ satisfying $\Bi{f_a,e_b}=\delta_{ab}$. A similar
observation holds for \eqref{eq:piGcoord}.}
\[
\tri_{G}= \frac{1}{12} \sum \Bi{e_a, [e_b, e_c]} e_a\wedge e_b\wedge e_c .
\]
A \textbf{quasi-Poisson $G$-manifold} \cite{AKM} consists of a
$G$-manifold $M$ together with a $G$-invariant bivector field
$\pi$ satisfying
\begin{equation}\label{eq:tripiG}
\tri_{\pi}= \rho_{M}(\tri_{G}),
\end{equation}
where $\rho_{M}: \mathfrak{g}\rmap \mathcal{X}(M)$ is the associated
infinitesimal action, and we keep the same notation for the
induced maps of exterior algebras. When $M$ is just a
$\mathfrak{g}$-manifold, we call the corresponding object a
\textbf{quasi-Poisson $\mathfrak{g}$-manifold}. The two notions
are related by the standard procedure of integration of
infinitesimal actions; in particular, they coincide if $M$ is
compact and $G$ is simply connected.

In analogy with ordinary or twisted Poisson manifolds, are
quasi-Poisson structures also associated with Lie algebroids? As
we now discuss, the answer is yes. Let us consider a more general
set-up: let $M$ be a $\mathfrak{g}$-manifold and let $\pi\in
\mathcal{X}^2(M)$ be an arbitrary bivector field. Motivated by
\cite[Sec.~9]{AKM}, we consider on $T^*M\oplus \mathfrak{g}$ the
``anchor'' map
\begin{equation}\label{eq:ranchor}
r: T^*M\oplus \mathfrak{g} \rmap TM, \;\; r(\alpha, v)=
\pi^{\sharp}(\alpha)+ \rho_{M}(v),
\end{equation}
combining the bivector field and the action. On sections of
$T^*M\oplus \mathfrak{g}$, we consider the bracket defined by
\begin{align}
& [(\alpha, 0), (\beta, 0)]= ([\alpha, \beta], 
\frac{1}{2}i_{\rho_{M}^{*}(\alpha\wedge\beta)}(\tri_{G})),\label{eq:brac1}\\
& [(0, v), (0, v')]= (0, [v, v']),\label{eq:brac2}\\
& [(0, v), (\alpha, 0)]= (\Lie_{\rho_M(v)}(\alpha),
0),\label{eq:brac3}
\end{align}
for all 1-forms $\alpha, \beta \in \Omega^1(M)$ and all $v, v'\in
\mathfrak{g}$ (thought of as constant sections in
$C^\infty(M,\mathfrak{g})$). As in Example
\ref{ex:transfalg}, the definition of the bracket on general
elements in $\Gamma(T^*M\oplus\mathfrak{g})=\Omega^1(M)\oplus
C^\infty(M,\mathfrak{g})$ is obtained from the Leibniz formula
\eqref{eq:leibniz}. With these definitions, we obtain a
quasi-Poisson analogue of Lemma \eqref{bivector}:

\begin{theorem}
Let $M$ be a $\mathfrak{g}$-manifold equipped with a bivector
field $\pi$. The following are equivalent:
\begin{enumerate}
\item[(i)] $(M,\pi)$ is a quasi-Poisson $\mathfrak{g}$-manifold;
 \item[(ii)] $r:\Omega^1(M)\oplus C^\infty(M,\mathfrak{g}) \to  \mathcal{X}(M)$ preserves
 brackets;
\item[(iii)] the bracket $[\cdot, \cdot]$  on $\Omega^1(M)\oplus
C^\infty(M,\mathfrak{g})$ satisfies the Jacobi identity;
 \item[(iv)]
$(T^*M\oplus\mathfrak{g}, r, [\cdot, \cdot])$ is a Lie algebroid.
\end{enumerate}
\end{theorem}

\begin{proof}
Note that $r$ preserves the bracket \eqref{eq:brac2}, since
$\rho_M$ is an action. From the identity \eqref{eq:ident1} in
Lemma \ref{bivector}, it follows that $r$ preserves the bracket of
type \eqref{eq:brac1} if and only if $\tri_{\pi}=
\rho_M(\tri_{G})$. On the other hand, $r$ preserves the bracket of
 type \eqref{eq:brac3} if and only if
$\pi^{\sharp}\Lie_{\rho_M(v)}(\xi)=
\Lie_{\rho_M(v)}\pi^{\sharp}(\xi)$, which is equivalent to the
$\mathfrak{g}$-invariance of $\pi$. This shows that (i) and (ii)
are equivalent.

Let us prove that (i) implies (iii); from the proof, the converse
will be clear. Assuming (i), we must show that $[\cdot,\cdot]$ on
$\Omega^1(M)\oplus C^\infty(M,\mathfrak{g})$ satisfies the Jacobi
identity. On elements of type $(0, v)$, this reduces to the Jacobi
identity for $\mathfrak{g}$ (or, alternatively, for
$\mathfrak{g}\ltimes M$). On elements $(0, v)$, $(0, w)$ and
$(\alpha, 0)$, the Jacobi identity of $[\cdot,\cdot]$ reduces to
the fact that $\rho_{M}$ is an action. Computing the
``jacobiator'' for elements of type $(0, v), (\alpha, 0), (\beta,
0)$, we see that the first component is
\begin{equation}\label{eq:jac1}
[\Lie_{\rho_M(v)}(\alpha), \beta]+ [\alpha,
\Lie_{\rho_M(v)}(\beta)]- \Lie_{\rho_M(v)}([\alpha, \beta]).
\end{equation}
Using the Leibniz identity, we see that the
$C^{\infty}(M)$-linearity of \eqref{eq:jac1} with respect to
$\beta$ is equivalent to $\pi^{\sharp}\Lie_{\rho(v)}(\beta)=
\Lie_{\rho(v)}\pi^{\sharp}(\beta)$, i.e., to the
$\mathfrak{g}$-invariance of $\pi$. Hence, if $\pi$ is invariant,
\eqref{eq:jac1} is $C^{\infty}(M)$-linear on $\alpha$ and $\beta$,
and then one can check that it is zero by looking at the
particular case when $\alpha$ and $\beta$ are exact. The second
component of the jacobiator of $(0, v), (\alpha, 0), (\beta, 0)$
can be computed similarly.

To complete the proof that (i) implies (iii), we must deal with
the Jacobi identity for elements of type $(\alpha, 0)$, $(\beta,
0)$, $(\gamma, 0)$. To this end, we first need to find the
expression for the bracket between elements of type $(0,
\tilde{v})$ and $(\alpha, 0)$, with $\tilde{v}\in C^{\infty}(M,
\mathfrak{g})$ not necessarily constant: pairing $d\tilde{v}\in
\Omega^1(M; \mathfrak{g})$ with an element $\mu\in C^{\infty}(M,
\mathfrak{g}^*)$ gives us a 1-form on $M$, denoted by
$A_{\tilde{v}}(\mu)$, satisfying the following two properties:
\[
A_{f\tilde{v}}(\mu)= fA_{\tilde{v}}(\mu)+ \mu(\tilde{v}) df,
\;\;\mbox{ and }\;\;\; A_{\tilde{v}}(f\mu)= f A_{\tilde{v}}(\mu),
\]
for $f\in C^\infty(M)$. We claim that
\begin{equation}
\label{Av} [(0, \tilde{v}), (\alpha, 0)]=
(\Lie_{\rho_M(\tilde{v})}(\alpha)-
A_{\tilde{v}}(\rho_{M}^{*}(\alpha)),
-\Lie_{\pi^{\sharp}(\alpha)}(\tilde{v})) .
\end{equation}
To see that, note that  \eqref{Av} holds when $\tilde{v}$ is
constant, and the difference between the left and right hand sides
is $C^{\infty}(M)$-linear in $\tilde{v}$. We remark that
\begin{equation}\label{formula-A-v}
A_{i_{\mu\wedge \mu'}(\tri_G)}(\mu'') + c.p. = 2d(\tri_{G}(\mu,
\mu', \mu'')).
\end{equation}
Again, it is easy to check this identity when $\mu$, $\mu'$ and
$\mu''$ are constant, so \eqref{formula-A-v} follows from
$C^{\infty}(M)$-linearity. Also, denoting $\tri_M:=
\rho_{M}(\tri_{G})$, a direct computation shows that
\[
\rho_{M}(i_{\rho_{M}^{*} (\alpha\wedge\beta)}(\tri_G))=
i_{\alpha\wedge\beta}(\tri_{M}).
\]
We now turn to the computation of the jacobiator of the elements
$(\alpha,0), (\beta, 0)$ and $(\gamma,0)$, that we denote by
$\Jac(\alpha, \beta, \gamma)$. For the first component of
$\Jac(\alpha,\beta,\gamma)$, we obtain
\begin{equation}\label{eq:jac31}
([\alpha,[\beta, \gamma]] + c.p.) -
\frac{1}{2}(\Lie_{i_{\alpha\wedge\beta}(\tri_M)}(\gamma) + c.p.) +
\frac{1}{2}(A_{i_{\alpha\wedge\beta}(\tri_M)}(\rho_M^*(\gamma))
+ c.p.).
\end{equation}
Combining the second identity of Lemma \ref{bivector} with
(\ref{formula-A-v}), we get that \eqref{eq:jac31} equals
\[
d((\rho_{M}(\tri_G)- \tri_{\pi})(\alpha, \beta, \gamma)),
\]
which vanishes by the condition $\tri_{\pi}= \rho_{M}(\tri_G)$. So
we are left with proving that the second component of
$\Jac(\alpha,\beta, \gamma)$ vanishes, which amounts to showing
that
\begin{equation}\label{eq:jac2}
i_{\rho_{M}^{*}([\alpha, \beta]\wedge \gamma)}(\tri_G) + c.p. =
\Lie_{\pi^{\sharp}(\gamma)}
i_{\rho_{M}^{*}(\alpha\wedge\beta)}(\tri_{G}) + c.p. .
\end{equation}
In order to do that, consider the operators $i_{\rho_M^*([\alpha,
\beta])}$ and $\Lie_{\pi^{\sharp}(\alpha)}i_{\rho_M^*(\beta)}-
\Lie_{\pi^{\sharp}(\beta)}i_{\rho_M^*(\alpha)}$  acting on
$C^{\infty}(M,\Lambda \mathfrak{g})$, for $\alpha,\beta \in
\Omega^1(M)$.
\begin{claim}
On $\Lambda \mathfrak{g}$, seen as constant functions in
$C^{\infty}(M,\Lambda \mathfrak{g})$, we have
\begin{equation}\label{i-L-formula}
i_{\rho_M^*([\alpha, \beta])} =
\Lie_{\pi^{\sharp}(\alpha)}i_{\rho_M^*(\beta)}-
\Lie_{\pi^{\sharp}(\beta)}i_{\rho_M^*(\alpha)}.
 \end{equation}
\end{claim}
\begin{proof}
Both operators are derivations of degree -1 on
$\Lambda\mathfrak{g}$, hence it suffices to show
\eqref{i-L-formula} for elements $v\in \mathfrak{g}$. As we now
check, this follows from the definition of the bracket induced by
$\pi$ and the invariance of $\pi$: on one hand,
\begin{equation}\label{eq:aux1}
i_{\rho_M^*([\alpha, \beta])}(v)= [\alpha, \beta](\rho_M(v)) =
i_{\rho_M(v)}\Lie_{\pi^{\sharp}(\alpha)}(\beta) -
i_{\rho_M(v)}\Lie_{\pi^{\sharp}(\beta)}(\alpha)-i_{\rho_M(v)}d\pi(\alpha,\beta).
\end{equation}
Using that $i_{[X,Y]}=\Lie_Xi_Y- i_Y\Lie_X$ for vector fields
$X,Y$, we have
\begin{eqnarray}\label{eq:aux2}
i_{\rho_M(v)}\Lie_{\pi^{\sharp}(\alpha)}(\beta)&=&
\Lie_{\pi^{\sharp}(\alpha)}(\beta(\rho_M(v)))-\beta([\pi^{\sharp}(\alpha),\rho_M(v)])\nonumber\\
&=&\Lie_{\pi^{\sharp}(\alpha)}i_{\rho_M^*(\beta)}(v)-\pi(\Lie_{\rho_M(v)}(\alpha),\beta),
\end{eqnarray}
where the last equality follows from the $\mathfrak{g}$-invariance
of $\pi$. Using the identity \eqref{eq:aux2} (and its analogue for
$\alpha$ and $\beta$ interchanged) in \eqref{eq:aux1},
\eqref{i-L-formula} follows.
\end{proof}

Using the claim, we see that
\begin{equation}\label{eq:aux4}
i_{\rho_M^*([\alpha, \beta]\wedge \gamma)} + c.p. =
i_{\rho_M^*(\gamma)}i_{\rho_M^*([\alpha, \beta])} + c.p.=
(i_{\rho_M^*(\gamma)}\Lie_{\pi^{\sharp}(\alpha)}i_{\rho_M^*(\beta)}-
i_{\rho_M^*(\gamma)}\Lie_{\pi^{\sharp}(\beta)}i_{\rho_M^*(\alpha)})
+ c.p.
\end{equation}
when restricted to constant elements in $C^{\infty}(M,
\Lambda\mathfrak{g})$ . On the other hand, it follows from
(\ref{i-L-formula}) that, on $\Lambda \mathfrak{g}$, we can write
\begin{equation}\label{eq:aux3}
i_{\rho_M^*([\alpha, \beta])}=
[\Lie_{\pi^{\sharp}(\alpha)},i_{\rho_M^*(\beta)}]-
[\Lie_{\pi^{\sharp}(\beta)},i_{\rho_M^*(\alpha)}]
\end{equation}
since the Lie derivatives are zero on constant functions. But both
sides of \eqref{eq:aux3} are $C^\infty(M)$-linear, so this
equality is valid for all $C^\infty(M,\mathfrak{g})$. So we can
write
\begin{eqnarray*}
i_{\rho_M^*([\alpha, \beta]\wedge \gamma)} + c.p.& =& -
i_{\rho_M^*([\alpha, \beta])}i_{\rho_M^*(\gamma)} + c.p.\\
& = & - ([\Lie_{\pi^{\sharp}(\alpha)},i_{\rho_M^*(\beta)}]-
[\Lie_{\pi^{\sharp}(\beta)},i_{\rho_M^*(\alpha)}])i_{\rho_M^*(\gamma)}
+ c.p.,
\end{eqnarray*}
from where we deduce that
\[
i_{\rho_M^*([\alpha, \beta]\wedge \gamma)} + c.p. = 2(
\Lie_{\pi^{\sharp}(\alpha)}i_{\rho_M^*(\beta\wedge\gamma)} +
c.p.)-
(i_{\rho_M^*(\gamma)}\Lie_{\pi^{\sharp}(\alpha)}i_{\rho_M^*(\beta)}-
i_{\rho_M^*(\beta)}\Lie_{\pi^{\sharp}(\alpha)}i_{\rho_M^*(\gamma)}
+ c.p.).
\]
On constant functions, we can use \eqref{eq:aux4} to conclude that
\[
 (i_{\rho_M^*([\alpha, \beta]\wedge \gamma)}+c.p.)= 2 (\Lie_{\pi^{\sharp}(\alpha)}i_{\rho_M^*(\beta\wedge\gamma)} + c.p.)-
(i_{\rho_M^*([\alpha, \beta]\wedge \gamma)}+c.p.),
\]
i.e., $i_{\rho_M^*([\alpha, \beta]\wedge \gamma)} + c.p. =
\Lie_{\pi^{\sharp}(\alpha)}i_{\rho_M^*(\beta\wedge\gamma)}$.
Evaluating this identity at $\tri_{G}$ proves \eqref{eq:jac2},
showing that (i) implies (iii). Looking back into the proof, one
can check that the same formulas show the converse, so that (i)
and (iii) are equivalent.

Since (iii) and the Leibniz identity for $[\cdot,\cdot]$ are together equivalent to (iv), it follows
that (i) -- (iv) are equivalent to each other.
\end{proof}

\begin{corollary}
\label{the-foliation}
If $(M,\pi)$ is a quasi-Poisson $\mathfrak{g}$-manifold, then the
generalized distribution
$$
\pi^{\sharp}(\alpha)+\rho_M(v) \subseteq TM, \;\mbox{ for }\;
\alpha \in T^*M,\;\; v \in \mathfrak{g},
$$
is integrable.
\end{corollary}
This result shows that the singular distribution discussed in
\cite[Thm.~9.2]{AKM} in the context of \textit{hamiltonian}
quasi-Poisson manifolds is integrable even without the presence of
a moment map (and without the positivity of $\Bi{\cdot,\cdot}$). 
As in the case of ordinary Poisson manifolds, we
call a quasi-Poisson manifold \textbf{nondegenerate} if its
associated Lie algebroid is transitive (i.e., its anchor map is onto).

\begin{example}\textit{(Quasi-Poisson structures on Lie groups)}\label{ex:quasip}

Let $G$ be a Lie group with Lie algebra $\mathfrak{g}$, which we
assume to be equipped with an invariant nondegenerate quadratic form
$\Bi{\cdot,\cdot}$. We consider $G$ acting on itself by
conjugation. As shown in \cite[Sec.~3]{AKM}, the bivector field $\pi_G$, defined on
left invariant 1-forms by
\begin{equation}\label{eq:quasip}
\pi_G(dl_{g^{-1}}^*(\mu),dl_{g^{-1}}^*(\nu)):=\frac{1}{2}\Bi{(\Ad_{g^{-1}}-\Ad_g)(\mu^\vee),\nu^\vee},
\end{equation}
where $l_g$ denotes left multiplication by $g\in G$, $\mu, \nu \in \mathfrak{g}^*$, and
$\mu^\vee$ is the element in $\mathfrak{g}$ dual to $\mu$ via $\Bi{\cdot,\cdot}$,
makes $G$ into a quasi-Poisson $G$-manifold. If $\Bi{\cdot,\cdot}$ is a metric, then we can write
\begin{equation}\label{eq:piGcoord}
\pi_G= \frac{1}{2}\sum e_a^l\wedge e_a^r,
\end{equation}
where $e_a$ is an orthonormal basis of $\mathfrak{g}$ and $e_a^r$
(resp. $e_a^l$) are the corresponding right (resp. left)
translations.

In this example, the image of $\pi_G^{\sharp}$ is tangent to the $G$-orbits,
so the leaves of the corresponding foliation are the conjugacy
classes. The formula for the Lie algebroid bracket on $T^*G\oplus \mathfrak{g}$
bears close resemblance with the one for the bracket in the ``double''
of the Lie quasi-bialgebra of $G$, as in  \cite{BK-S}.
We will discuss this connection 
in a separate work.
\end{example}

\section{Moment maps in Dirac geometry: the infinitesimal
picture}\label{sec:infinit}

\subsection{Dirac manifolds}
Let $\phi$ be a closed 3-form on a manifold $M$. A
\textbf{$\phi$-twisted Dirac structure} on $M$ \cite{SeWe01} is a
subbundle $L\subset E= TM \oplus T^*M$ satisfying the following
two conditions:
\begin{itemize}
\item[1.] $L$ is maximal isotropic with respect to the symmetric
pairing $\SP{\cdot,\cdot}_{\st{+}}: \Gamma(E) \times \Gamma(E) \to
C^\infty(M)$,
\begin{equation}\label{eq:brk1}
\SP{(X,\alpha),(Y,\beta)}_{\st{+}}:=\beta(X) + \alpha(Y);
\end{equation}
\item[2.] The space of sections $\Gamma(L)$ is closed under the
bracket $\Cour{\cdot,\cdot}_{\phi}:\Gamma(E)\times \Gamma(E) \to
\Gamma(E)$,
\begin{equation}\label{eq:brk2}
\Cour{(X,\alpha),(Y,\beta)}_{\phi}:= ([X,Y],
\Lie_X\beta-i_Yd\alpha +i_{X\wedge Y}\phi).
\end{equation}
\end{itemize}
Since the pairing \eqref{eq:brk1} has zero signature, condition 1.
is equivalent to requiring that $L$ has rank equal to $\dim(M)$
and that $\SP{\cdot,\cdot}_{\st{+}}|_L=0$. The bracket \eqref{eq:brk2} is
the \textbf{$\phi$-twisted Courant bracket} considered in
\cite{SeWe01}. When $\phi$=0, this bracket is a non-skew-symmetric
version of Courant's original bracket introduced in \cite{Cou90}.

Twisted Dirac structures are
always associated with Lie algebroids. Indeed,
the restriction of the Courant
bracket $\Cour{\cdot,\cdot}_{\phi}$ to a Dirac subbundle $L\subset
TM \oplus T^*M$ defines a Lie algebra bracket on the space of sections $\Gamma(L)$,
making $L \to M$ into  a Lie algebroid with anchor
$$
\rho=\pr_1|_L:L\to TM,
$$
where $\pr_1$ is the first projection. The orbits of this
algebroid are also called the \textbf{leaves} of $L$.

\begin{example} (\textit{Twisted Poisson structures})

If $\pi$ is a bivector field on $M$, then
$$
L_{\pi}:=\gra({\pi}^{\sharp}) \subset TM \oplus T^*M
$$
satisfies
condition 1., and $L_{\pi}$ is a $\phi$-twisted Dirac structure if
and only if $\pi$ is a $\phi$-twisted Poisson structure in the sense of
Example \ref{ex:twispoiss}. In this case, the second projection
$$
\pr_2|_L:L \to T^*M
$$
establishes an isomorphism of Lie algebroids, where $T^*M$ is
equipped with the Lie algebroid structure described in Example
\ref{ex:twispoiss}. Setting $\phi=0$, we obtain a one-to-one
correspondence between ordinary Poisson structures on $M$ and
Dirac structures satisfying the extra condition $L \cap TM
=\{0\}$.
\end{example}
\begin{example} (\textit{Twisted presymplectic forms})

Similarly, the graph associated with a 2-form $\omega \in
\Omega^2(M)$, $L_{\omega}=\gra({\omega}^{\sharp})$, is a
$\phi$-twisted Dirac structure if and only if
$$
d\omega + \phi =0,
$$
and we refer to $\omega$ as a \textbf{$\phi$-twisted
presymplectic form}. In this case, setting $\phi=0$, we have an identification
of closed 2-forms
on $M$ with Dirac structures
satisfying $L\cap T^*M=\{0\}$.
\end{example}

In general, the leaves of a $\phi$-twisted Dirac structure $L$
carry twisted presymplectic forms defined as follows: at each $x
\in M$, we define a skew symmetric bilinear form $\theta_x$ on
$\rho(L)_x=\pr_1(L)_x$ by
\begin{equation}\label{eq:2form}
\theta_x(X_1,X_2)=\alpha(X_2),
\end{equation}
where $\alpha$ is any element in $T_x^*M$ satisfying $(X_1,\alpha)
\in L_x$. The fact that $L$ is maximal isotropic with respect
to \eqref{eq:brk1} guarantees that \eqref{eq:2form} is independent of the choice of $\alpha$, and
these forms fit together into a smooth
leafwise 2-form $\theta$.  Using that
$\Gamma(L)$ is closed with respect
to \eqref{eq:brk2}, one can show that, on each leaf
$\iota:\mathcal{O}\hookrightarrow M$, the 2-form $\theta$ satisfies
$$
d\theta + \iota^*\phi=0.
$$
At each point $x\in M$, the kernel of $\theta$ coincides with
$L_x\cap T_xM$, which shows that the leafwise presymplectic forms
are nondegenerate if and only if $L$ comes from a $\phi$-twisted
Poisson structure. We will denote the distribution $L\cap TM$ on $M$
by $\ker(L)$.

Since Dirac structures are always associated with Lie algebroids,
it is natural to consider how to obtain Dirac structures from
them. The following is a useful construction, see \cite{BCWZ}: for
a Lie algebroid $A$ over $M$ with anchor $\rho: A\rmap TM$, we
define a \textbf{$\phi$-IM form of $A$} to be  any bundle map\footnote{These bundle maps are
infinitesimal versions of multiplicative 2-forms on groupoids, see
\cite{BCWZ}; the terminology ``IM'' stands for ``infinitesimal
multiplicative''.}
\[
\sigma: A\rmap T^*M
\]
satisfying the following properties:
\begin{eqnarray}
\SP{\sigma(\xi),\rho(\xi')} & = & -\SP{\sigma(\xi'),\rho(\xi)};\label{eq:skew}\\
\sigma([\xi, \xi']) & = & \mathcal{L}_{\xi}(\sigma(\xi'))-
                                      \mathcal{L}_{\xi'}(\sigma(\xi))+
                                      d\SP{\sigma(\xi),\rho(\xi')} +
                                      i_{\rho(\xi)\wedge
                                      \rho(\xi')}(\phi),\label{eq:Dirac}
\end{eqnarray}
for $\xi,\xi' \in \Gamma(A)$ (here $\SP{\cdot,\cdot}$ denotes the usual pairing between a 
vector space and its dual). Let $L_{\sigma}\subset TM\oplus
T^*M$ be the image of the map $(\rho, \sigma):A \rmap TM\oplus T^*M$.
Then the following is immediate.
\begin{lemma}
\label{multiplicat} If $\sigma$ is a $\phi$-IM form of $A$ and
$\mathrm{rank}(L_{\sigma})= \dim(M)$, then $L_{\sigma}$ is a
$\phi$-twisted Dirac structure on $M$.
\end{lemma}
Of course, any Dirac structure can be realized as the image of an IM form by
taking $A= L$, viewed as an algebroid with $\rho= \pr_1|L$, and
$\sigma= \pr_2|L$.

The following is a key example.

\begin{example}(\textit{Cartan-Dirac structures on Lie groups})\label{ex:CD}

Cartan-Dirac structures on Lie groups play a role in Dirac
geometry analogous to the one played by Lie-Poisson structures
(Example \ref{Lie-Poisson}) in Poisson geometry. Just as
Lie-Poisson structures on the dual of Lie algebras are completely
determined by the Kostant-Kirillov-Souriau (KKS) symplectic forms
along coadjoint orbits, Cartan-Dirac structures on Lie groups
``assemble'' certain 2-forms defined on conjugacy classes defined
as follows.

Let $G$ be a Lie group with Lie algebra $\mathfrak{g}$, and let
$\Bi{\cdot,\cdot}$ be a bi-invariant nondegenerate quadratic form, which we use to
identify $TG$ and $T^*G$. For $v \in \mathfrak{g}$, let
$v_G=v_r-v_l$ be the infinitesimal generator of the action of $G$
on itself by conjugation. We define, on each conjugacy class
$\iota: \mathcal{C} \hookrightarrow G$, a 2-form $\theta$ by
\begin{equation}\label{eq:GHJW}
\theta_g(u_G,v_G):=
(\frac{1}{2}(\Ad_{g}-\Ad_{g^{-1}})u,v)_{\mathfrak{g}}, \;\;\; g
\in \mathcal{C}.
\end{equation}
Direct computations show that $d\theta - \iota^*\phi^{G}=0$, where
$\phi^{G}$ is the bi-invariant Cartan 3-form on $G$, and that
$\theta_g$ is nondegenerate at a point $g$ if and only if $(\Ad_g
+ 1)$ is invertible. The 2-forms \eqref{eq:GHJW} appear in \cite{GHJW} in the
study of symplectic structures of moduli spaces.

Since these 2-forms are not symplectic, but \textit{twisted
presymplectic}, they should correspond to a $-\phiG$-twisted Dirac
structure $L_G$ on $G$ rather than a Poisson structure. A simple
computation shows that
\begin{equation}\label{eq:cartandirac}
L_G=\{ (v_r-v_l, \frac{1}{2}(v_r+ v_l)),\;\; v\in \mathfrak{g}\}
\subset TG\oplus TG.
\end{equation}
(Recall that
we are identifying $TG$ with $T^*G$ via $\Bi{\cdot,\cdot}$.)
We call $L_G$ the \textbf{Cartan-Dirac structure} on $G$
associated with $\Bi{\cdot,\cdot}$ \cite{SeWe01,BCWZ}.

Note that $\rho(v) = v_r-v_l$ is the anchor of the action Lie
algebroid (Example \ref{ex:transfalg}) $\mathfrak{g}\ltimes G$
with respect to the action by conjugation, and the map
\begin{equation}
\label{infinitesimal-Dirac} \sigma: G\times \mathfrak{g}\rmap TG,
\;\; \sigma(v)=  \frac{1}{2}(v_r+ v_l)
\end{equation}
satisfies the conditions of Lemma \ref{multiplicat}. So
$\sigma$ is a $-\phi^{G}$-IM form of  $\mathfrak{g}\ltimes G$, and
the Cartan-Dirac structure $L_G$ arises as the image of
$(\rho,\sigma)$. In this case, $(\rho,\sigma)$ actually
establishes an isomorphism between $\mathfrak{g}\ltimes G$ and
$L_G$. (Note the analogy with Example \ref{Lie-Poisson}, which
shows that Lie algebroids of Lie-Poisson structures are isomorphic
to action Lie algebroids for the coadjoint action!)
\end{example}

Let us finally recall an important operation involving Dirac structures:
if $L$ is a $\phi$-twisted Dirac structure on $M$ and $B \in \Omega^2(M)$,
then
\begin{equation}
\tau_B(L):=\{(X,\alpha + {B}^{\sharp}(X))\;|\; (X,\alpha) \in L\}
\end{equation}
defines a $(\phi-dB)$-twisted Dirac structure on $M$
\cite{SeWe01}. The operation $\tau_B$ is called a \textbf{gauge
transformation associated with $B$}, and it has the effect of
modifying $L$ by adding the pull-back of $B$ to the presymplectic
form on each leaf.

\subsection{Dirac maps}\label{subsec:diracmaps}

Since Dirac structures generalize both Poisson and
presymplectic structures, we have two possible definitions of Dirac maps,
see \cite{BuRa02}.

Let $(M,\LM)$ and $(N,\LN)$ be twisted Dirac manifolds.
A smooth map $f: N \to M$ is a \textbf{forward Dirac map},
or \textbf{f-Dirac} in short, if $\LN$ and $\LM$ are related as
follows:
\begin{equation}\label{eq:fdirac}
\LM=\{(df(Y),\alpha)\;|\; Y\in TN,\, \alpha \in T^*M
\mbox{ and } (Y,df^*(\alpha)) \in \LN \}.
\end{equation}
If $\LM$ and $\LN$ are associated with twisted Poisson
structures, then an f-Dirac map is equivalent to a Poisson map.
The terminology ``forward'' is due to the  fact that, at each
point, \eqref{eq:fdirac} extends the usual notion of
``push-forward'' of a linear bivector. For this reason, we may
write
$$
\LM=f_*\LN
$$
instead of \eqref{eq:fdirac}, in analogy with the notation for ``$f$-related''
bivector fields on a manifold.

Similarly, $f: N \to M$ is a \textbf{backward Dirac map}, or
simply \textbf{b-Dirac}, if
\begin{equation}\label{eq:bdirac}
\LN=\{(Y,df^*\alpha)\;|\; Y\in TN,\, \alpha \in T^*M
\mbox{ and } (df(Y),\alpha) \in \LM\}.
\end{equation}
If $\LM$ and $\LN$ are associated with twisted
presymplectic structures $\omega_{\st{M}}$ and $\omega_{\st{N}}$,
then a b-Dirac map is just a map satisfying
$f^*\omega_{\st{M}} = \omega_{\st{N}}$. As before, we will
write
$$
\LN=f^*\LM
$$
to denote that $f$ is a b-Dirac map. Note that
$f^*\LM$ is always a well-defined, though not necessarily
smooth, subbundle of $TN$, in contrast with $f_*\LN$, which
may not be well-defined at all. In fact, $f^*\LM$ defines a
Dirac structure on $N$ provided it is smooth, which is the case,
e.g., when $f$ is a submersion. However, as illustrated in
the next example, $f^*\LM$ may define a Dirac structure even
when $f$ is not a submersion.

\begin{example}(\textit{Inclusion of presymplectic leaves})\label{ex:incl}

Let $L$ be a twisted Dirac structure on $M$, and consider a
presymplectic leaf $\mathcal{O}$, equipped with
Dirac structure $L_{\theta}$ associated with the twisted presymplectic form
$\theta$. Denoting by $\iota: \mathcal{O} \hookrightarrow M$ the inclusion map,
it follows from the definition of $\theta$ that
\begin{equation}\label{eq:bincl}
L_{\theta}=\{(X,i_X\theta)\;|\; X \in T\mathcal{O}   \} =
\{(X,(d\iota)^*\alpha)\;|\; (d\iota(X),\alpha) \in L\} = \iota^*L.
\end{equation}
So $\iota: (\mathcal{O},L_{\theta}) \hookrightarrow (M,L)$ is a b-Dirac map.
On the other hand, at each point of $M$, we have
$$
\iota_*L_{\theta}= \{(d\iota(X),\alpha)\;|\;
(X,(d\iota)^*\alpha) \in L_{\theta}\}.
$$
By the second equality in \eqref{eq:bincl}, it follows that $\iota_*L_{\theta} \subseteq L$,
but since they have the same dimension, we get
\begin{equation}\label{fincl}
\iota_*L_{\theta}=L,
\end{equation}
so $\iota$ is f-Dirac as well.
\end{example}
Note that the fact that the inclusion of presymplectic leaves into a
Dirac manifold is an f-Dirac map
is a direct generalization of
the fact that the inclusion of symplectic leaves into a Poisson manifold is a Poisson map.
As a simple consequence, we have
\begin{corollary}\label{cor:incl}
Let $(N,\LN)$ and $(M,\LM)$ be twisted Dirac manifolds. A map $J:N\to M$
is f-Dirac if and only if its restriction to each presymplectic leaf of $N$
is f-Dirac.
\end{corollary}
We remark that Example \ref{ex:incl} is very special in that the
inclusion map of presymplectic leaves is both forward and backward
Dirac (see also Remark \ref{rem:bandf}). In general, f-Dirac maps need not be b-Dirac, nor the other
way around.

\subsection{Poisson maps as infinitesimal hamiltonian actions}

The usual notion of Lie algebra action can be extended to the realm of Lie algebroids,
the main difference being that algebroids, rather than acting on manifolds,
act on maps from manifolds into their base \cite{MaHi}:
An \textbf{action of a Lie algebroid $A \to M$ on a map $J: N \to M$} consists of
a Lie algebra homomorphism $\rho_N: \Gamma(A) \to \mathcal{X}(N)$ satisfying
\begin{equation}\label{eq:equiva}
dJ\circ \rho_N (\xi)  = \rho(\xi), \;\; \mbox{ for all } \, \xi \in \Gamma(A),
\end{equation}
and such that, for $f\in C^\infty(M)$ and $\xi \in \Gamma(A)$,
$\rho_N(f\xi)= J^*f \rho_N(\xi)$ (i.e., the induced map
$\Gamma(J^*A)\to \mathcal{X}(N)$ comes from a vector bundle
morphism $J^*A \to TN$, where $J^*A=A\times_M N$ is the pull-back
of the vector bundle $A$ by $J$).

\begin{example}(\textit{Actions of transformations Lie algebroids})\label{ex:actionstransf}

Consider an infinitesimal action $\rho$ of $\mathfrak{g}$ on a manifold $M$.
Then an action $\rho_N$ of the transformation Lie algebroid $A=\mathfrak{g}\ltimes M$ on a map $J:N\to M$
is equivalent to an infinitesimal action $\overline{\rho_N}$ of $\mathfrak{g}$ on $N$ for
which $J$ is $\mathfrak{g}$-equivariant. Indeed, $\rho_N$ and $\overline{\rho_N}$ are related by
the formula
\begin{equation}
\rho_N(v)_y= \overline{\rho_N}(v(J(y)))_y,\;\; \mbox{ where } v\in C^\infty(M,\mathfrak{g}),\, y \in N,
\end{equation}
and the $\mathfrak{g}$-equivariance of $J$ corresponds to \eqref{eq:equiva}.
\end{example}

In Poisson geometry, Poisson maps are \textit{always} associated with Lie algebroid actions:
If $(Q,\piQ)$ and $(P,\piP)$ are Poisson manifolds, then any Poisson map $J:Q \to P$
induces a Lie algebroid action of $T^*P$ on $Q$ by
\begin{equation}
\label{eq:actionpoissonmap}
\Omega^1(P) \longrightarrow \mathcal{X}(Q),\;\;\; \alpha \mapsto
{\piQ}^{\sharp}(J^*\alpha).
\end{equation}
When the target $P$ is the dual of a Lie algebra, we recover a
familiar example:

\begin{example}(\textit{Infinitesimal hamiltonian actions})\label{ex:infinitham}

Consider $\mathfrak{g}^*$ equipped with
its Lie-Poisson structure. As remarked in Example \ref{Lie-Poisson},
the Lie algebroid structure on $T^*\mathfrak{g}^*$ induced by
\eqref{eq:liepoiss} is that of a transformation Lie algebroid $\mathfrak{g}\ltimes \mathfrak{g}^*$
with respect to the coadjoint action. If $J:Q \to \mathfrak{g}^*$ is a Poisson map, then it induces
an action of $T^*\mathfrak{g}^*$ via \eqref{eq:actionpoissonmap}, which, by Example \ref{ex:actionstransf},
is equivalent to an ordinary $\mathfrak{g}$-action on $Q$ for which $J$ is equivariant.
A simple computation shows that the $\mathfrak{g}$-action  arising in this way is just a hamiltonian
action in the usual sense, making $Q$ into a hamiltonian Poisson $\mathfrak{g}$-manifold
having $J$ as a momentum map.
\end{example}

Recall that a Poisson map $J:Q \to P$ is called a \textbf{symplectic realization} if $Q$ is
symplectic. The following is an immediate consequence:

\begin{proposition}\label{prop:1-1}
There is a one-to-one correspondence between Poisson maps into
$\mathfrak{g}^*$ and hamiltonian Poisson $\mathfrak{g}$-manifolds, and this correspondence restricts
to a one-to-one correspondence between symplectic realizations of $\mathfrak{g}^*$
and hamiltonian symplectic $\mathfrak{g}$-manifolds.
\end{proposition}

\begin{remark}\label{rem:poisslie1}
An analogue of Prop.~\ref{prop:1-1} holds more generally in the context
of Poisson-Lie groups \cite{Luthesis,LuWe}. Let $(G,\pi)$
be a simply-connected Poisson-Lie group, and let $G^*$ be its dual. The Lie algebroid
structure on $T^*G^*\cong G^* \times \mathfrak{g}$ induced from the dual Poisson structure is a transformation
Lie algebroid, now associated with the infinitesimal
dressing action of $\mathfrak{g}$ on $G^*$. For a Poisson map $J:Q \to G^*$,
the general Lie algebroid action described by \eqref{eq:actionpoissonmap} reduces to a
Poisson $\mathfrak{g}$-action for which $J$ is an equivariant momentum map in the sense of Lu \cite{Lu}.
\end{remark}

In order to extend this discussion to Dirac geometry, let us consider
$L_{\piP} = \gra({\piP}^{\sharp})$,
the associated Dirac structure on $(P,\piP)$. Using the Lie algebroid isomorphism
$T^*P \cong L_{\piP}$, we can rewrite the infinitesimal action \eqref{eq:actionpoissonmap} as
\begin{equation}\label{eq:actionLpi}
\Gamma(L_{\piP}) \to \mathcal{X}(Q),\;\; (X,\alpha) \mapsto Y,
\end{equation}
where $Y \in \mathcal{X}(Q)$ is uniquely determined by the
condition $(Y,J^*\alpha) \in L_{\piQ}$. Also note that $Y$ is
related to $X$ by $dJ(Y)=X$, since $J$ is a Poisson map. The
question of whether this procedure can be carried out for
f-Dirac maps leads us to the notion of Dirac
realization.

\subsection{Dirac realizations}
If $(N,\LN)$ and $(M,\LM)$ are twisted Dirac manifolds, then, by
definition, a smooth $J:N \to M$ is an f-Dirac map if and only if,
given $(X,\alpha) \in (\LM)_{J(y)}$, there exists a $Y \in T_yN$
with the property that
\begin{equation}\label{eq:fdircond}
(Y,dJ^*\alpha) \in (\LN)_y \;\; \mbox{ and } \; X=(dJ)_y(Y).
\end{equation}
It is natural to try to define an action of $\LM$ on $N$ just as in the case of
Poisson maps, see \eqref{eq:actionLpi}, except that \eqref{eq:fdircond} does \textit{not} determine
$Y$ uniquely in general. In fact, this is the case if and only if
the following extra ``nondegeneracy'' condition holds:
\begin{equation}\label{eq:nondegreal}
\ker(dJ)\cap \ker(\LN)=\{0\}.
\end{equation}
A similar argument as in \cite[Section~7.1]{BCWZ} shows that \eqref{eq:nondegreal} is equivalent to
$J:\ker(\LN) \to \ker(\LM)$ being an isomorphism.

\begin{definition}\label{def:diracreal}
A \textbf{Dirac realization} of a $\phi$-twisted Dirac manifold
$(M,\LM)$ is an f-Dirac map $J:(N,\LN)\to (M,\LM)$, where $\LN$ is
a $J^*\phi$-twisted Dirac structure on $N$, satisfying
\eqref{eq:nondegreal}.
\end{definition}

As a consequence of Definition \ref{def:diracreal}, we have
\begin{corollary}
Let $J:N\to M$ be a Dirac realization. Then the map $\Gamma(\LM)
\to \mathcal{X}(N)$, $(X,\alpha) \mapsto Y$, where $Y$ is
determined by the conditions in \eqref{eq:fdircond}, is a Lie
algebroid action.
\end{corollary}

Dirac realizations $J:N \to M$ for which $N$ is presymplectic
were studied in \cite[Sec.~7.1]{BCWZ} under the name of
\textbf{presymplectic realizations}.
As a result of Corollary \ref{cor:incl}, we have
\begin{corollary}\label{cor:restric}
A map $J:N\to M$ is a Dirac realization if and only if its restriction
to each presymplectic leaf of $N$ is a presymplectic realization.
\end{corollary}

Similarly to Poisson geometry, the connection between Dirac realizations
and ``hamiltonian actions'' is established by a suitable choice of
``target'' $M$. Following the analogy between Lie-Poisson structures and Cartan-Dirac structures, it is natural
to study the ``hamiltonian spaces'' associated with Dirac realizations of Cartan-Dirac structures.
The particular case of presymplectic realizations is discussed in \cite[Sec.~7.2]{BCWZ}:

\begin{example}(\textit{Presymplectic realizations of Cartan-Dirac structures})\label{ex:realCD}

Let $G$ be a Lie group with Lie algebra $\mathfrak{g}$, equipped with a bi-invariant nondegenerate
quadratic form $\Bi{\cdot,\cdot}$. Let $L_G$ be the associated Cartan-Dirac
structure on $G$. If $(M,\omega_M)$ is a twisted presymplectic
manifold, then the conditions for $J:M \to G$ being a
presymplectic realization can be expressed as follows:
\begin{enumerate}
\item $\omega_M$ is $\mathfrak{g}$-invariant and  satisfies
$d\omega_M= J^*\phi^{G}$;

\item at each $x \in M$, $\Ker(\omega_M)_x =  \{(\rho_{M})_x(v): v
\in \Ker(\Ad_{J(p)}+ 1)\}$;

\item the map $J$ satisfies the moment map condition
\begin{equation}
\label{mom-1}
\omega^{\sharp}\rho_{M}= J^*\sigma.
\end{equation}
\end{enumerate}
The invariance of $\omega_M$ in 1. is with respect to the
$\mathfrak{g}$-action $\rho_M$ induced by $J$ (recall that $L_G
\cong \mathfrak{g}\ltimes G$, see Example \ref{ex:CD}, so an $L_G$-action
defines an ordinary $\mathfrak{g}$-action), for which
$J$ is equivariant; in 3., $\sigma$ is the IM-form of the Cartan-Dirac
structure (\ref{infinitesimal-Dirac}),
\begin{equation}\label{eq:sigma}
\sigma: \mathfrak{g}\rmap T^*G,\;\;\; \sigma(v)= \frac{1}{2}(v_r+ v_l)^{\vee},
\end{equation}
where $v\rmap v^{\vee}$ denotes the isomorphism $TG\rmap T^*G$
induced by the quadratic form. The ``relative closedness'' of $\omega_M$
in 1. expresses that the associated Dirac structure is
$-J^*\phi^{G}$-twisted, while  condition 2. is the ``non-degeneracy''
condition (\ref{eq:nondegreal}) applied to this particular case; finally,
condition 3. follows from $J$ being an f-Dirac map.

Conditions 1., 2. and 3.  are exactly the defining axioms of a \textbf{quasi-hamiltonian $\mathfrak{g}$-space},
in the sense of \cite{AMM}, for which $J$ is the group-valued moment map. Conversely,
any group-valued moment map of a quasi-hamiltonian $\mathfrak{g}$-space is a presymplectic
realization of $(G,L_G)$.
\end{example}

We summarize Example \ref{ex:realCD} in the next result, analogous
to Proposition \ref{prop:1-1}, see \cite[Thm.~7.6]{BCWZ}.
\begin{theorem}
\label{quasi-hamiltonian}
There is a one-to-one correspondence between presymplectic realizations of
$G$ endowed with the Cartan-Dirac structure, and quasi-hamiltonian $\mathfrak{g}$-manifolds.
\end{theorem}

Combining Corollary \ref{cor:restric} with Theorem \ref{quasi-hamiltonian}, we conclude that general
Dirac realizations of Cartan-Dirac structures  must be ``foliated'' by quasi-hamiltonian
$\mathfrak{g}$-manifolds. Since hamiltonian quasi-Poisson manifolds, in the sense of \cite{AKM}, also
have this property \cite[Sec.~10]{AKM}, we are led to investigate the relationship between these objects.

\subsection{Dirac realizations and hamiltonian quasi-Poisson $\mathfrak{g}$-manifolds}

\subsubsection{The equivalence theorem}

For a quasi-Poisson $\mathfrak{g}$-manifold $(M, \pi)$,
a \textbf{momentum map} is a $\mathfrak{g}$-equivariant map $J: M\rmap G$ (with respect to the infinitesimal
action by conjugation on $G$) satisfying the condition \cite[Lem.~2.3]{AKM}
\begin{equation}
\label{mom-2}
\pi^{\sharp}J^*= \rho_{M}\sigma^{\vee},
\end{equation}
where
\begin{equation}\label{eq:sigmavee}
\sigma^{\vee}: T^*G\rmap \mathfrak{g},\;\;\;
\sigma^{\vee}(\xi_g) = \frac{1}{2}(dr_{g^{-1}}(\xi_{g}^{\vee})+ dl_{g^{-1}}(\xi_{g}^{\vee}))
\end{equation}
is the adjoint of $\sigma$ \eqref{eq:sigma} with respect to the form $\Bi{\cdot,\cdot}$, and $\rho_M:\mathfrak{g}
\to TM$ is the $\mathfrak{g}$-action. Here $l_g$ and $r_g$ denote the left and right
translations by $g$, respectively, and, as in Example \ref{ex:realCD}, the symbol $^{\vee}$ on elements
of $T^*G$ is used to denote the corresponding element in $TG$ via the identification induced by $\Bi{\cdot,\cdot}$
(and vice-versa).

The following is our main result in this section.

\begin{theorem}
\label{quasi-Poisson-hamiltonian}
There is a one-to-one correspondence between Dirac realizations of
$G$, endowed with the Cartan-Dirac structure, and hamiltonian quasi-Poisson $\mathfrak{g}$-manifolds.
\end{theorem}

Before proving Theorem \ref{quasi-Poisson-hamiltonian},
let us collect some useful formulas relating the maps  $\sigma$,
$\sigma^{\vee}$, $\rho: \mathfrak{g} \to TG$, $\rho(v) = v_r-v_l$, and, for symmetry,
the dual of $\rho$ with respect to $\Bi{\cdot,\cdot}$,
\begin{equation}\label{eq:rhovee}
\rho^{\vee}: TG\rmap \mathfrak{g},\;\;\;  \rho^{\vee}(V_g)= dr_{g^{-1}}(V_g)- dl_{g^{-1}}(V_g) .
\end{equation}
The following lemma follows from a  straightforward computation.
\begin{lemma} The following formulas hold true:
\begin{eqnarray}
4 \sigma^{\vee}\sigma+ \rho^{\vee}\rho & = & 4 \id_{\mathfrak{g}},\label{eq:id1}\\
4 \sigma\sigma^{\vee}+ (\rho\rho^{\vee})^{*} & = & 4 \id_{T^*G},\label{eq:id2}\\
\sigma^*\rho & = & - \rho^*\sigma,\label{eq:id3}\\
 \sigma\rho^{\vee} & = & - (\rho^{\vee})^{*}\sigma^*,\label{eq:id4}\\
\rho^{\vee}(\sigma^{\vee})^* & = &  - \sigma^{\vee}(\rho^{\vee})^*,\label{eq:id5}\\
\rho\sigma^\vee &=& -(\rho\sigma^\vee)^*.
\end{eqnarray}
\end{lemma}

Motivated by the equivalence between quasi-hamiltonian manifolds and
nondegenerate hamiltonian quasi-Poisson manifolds \cite[Thm.~10.3]{AKM}, which will
also follow from Theorem \ref{quasi-Poisson-hamiltonian}, it is
natural to combine the two moment-map conditions (\ref{mom-1}) and
(\ref{mom-2}). By applying $\pi^{\sharp}$ to (\ref{mom-1}) and
using (\ref{mom-2}), we obtain, in particular, the relation
$$
\rho_{M}\sigma^{\vee}\sigma= \pi^{\sharp}\omega^{\sharp}\rho_{M}.
$$
This suggests the importance of writing $\rho_{M}\sigma^{\vee}\sigma$  as the composition of some
operator $C:TM \to TM$ with $\rho_M$ in general. Using \eqref{eq:id1} and the equivariance of $J$,
written as $\rho= dJ\circ \rho_{M}$, it is easy to find an expression for $C$ (which already appears in
\cite[Lem.~10.2]{AKM}):

\begin{lemma}
For any manifold $M$ equipped with an infinitesimal action $\rho_M: \mathfrak{g}\rmap TM$, and any
$\mathfrak{g}$-equivariant map $J: M\rmap G$, the operator
\[
C= 1- \frac{1}{4}\rho_{M}\rho^{\vee}(dJ): TM\rmap TM ,
\]
and its dual $C^*:T^*M \rmap T^*M$, satisfy the formulas
\begin{equation}\label{eq:CC*}
\rho_{M}\sigma^{\vee}\sigma= C\rho_{M}, \;\mbox{ and }\;  J^*\sigma\sigma^{\vee}= C^*J^*.
\end{equation}
\end{lemma}

Theorem \ref{quasi-Poisson-hamiltonian} follows from the next two propositions, each one of them
describing explicitly one direction of the asserted  one-to-one correspondence.

\begin{proposition}
\label{prop-direct}
Let $M$ be a quasi-Poisson $\mathfrak{g}$-manifold, and let $A= T^*M\oplus \mathfrak{g}$ be its associated
Lie algebroid, with anchor $r$. Then any moment map $J: M \rmap G$ induces a $-J^*\phi^G$-IM form of $A$ by
\begin{equation}\label{eq:s}
 s: A\rmap T^*M, \;\;\;   s(\alpha, v)= C^*(\alpha)+ J^*\sigma(v),
\end{equation}
so that the image $L$ of the map $(r, s): A\rmap TM\oplus T^*M$ is a $-J^*\phi^{G}$-twisted Dirac structure on $M$, and
$J:(M,L) \rmap (G,L_G)$ is a Dirac realization of the Cartan-Dirac structure on  $G$.
\end{proposition}

This proposition also suggests the converse construction.
\begin{proposition}
\label{prop-inverse}
Let $J: (M,L) \rmap (G,L_G)$ be a Dirac realization of the Cartan-Dirac structure on $G$. Then
\begin{enumerate}
\item[(i)] for any $v\in \mathfrak{g}$ there is an unique vector $V\in TM$ satisfying
\begin{equation}
\label{construct-action}
dJ(V)= \rho(v),\;\; \mbox{ and }\;\;  (V, J^*\sigma(v)) \in L.
\end{equation}
\item[(ii)] for any $\alpha \in T^*M$, there is an unique vector $X\in TM$ satisfying
\begin{align}
& dJ(X)= - (\rho_M\sigma^{\vee})^*\alpha,\label{pi-cond1}\\
& (X, C^*(\alpha))\in L .\label{pi-cond2}
\end{align}
\end{enumerate}
Moreover, $v \mapsto \rho_M(v):= V$ defines a $\mathfrak{g}$-action on $M$, and
$\alpha \mapsto \pi^{\sharp}(\alpha):= X$ defines a quasi-Poisson tensor $\pi$ on $M$ so that
$(M, \pi)$ is a hamiltonian quasi-Poisson $\mathfrak{g}$-manifold with moment map $J$.
\end{proposition}
Note that the $\mathfrak{g}$-action defined by \eqref{construct-action} is just the one induced by
the infinitesimal $L_G=\mathfrak{g}\ltimes G$-action with ``moment'' $J:M\to G$.

It is simple to check that the constructions in Propositions \ref{prop-direct} and \ref{prop-inverse}
are inverses to one another. For example, if $L$ is obtained from $\pi$ as in Prop.~\ref{prop-direct},
then
$$
L=\{(\pi^{\sharp}(\alpha) + \rho_M(v), C^*(\alpha) + J^*\sigma(v))\,|\, \alpha \in T^*M, v\in
\mathfrak{g} \},
$$
and it is clear that, given $\alpha \in T^*M$, $X = \pi^{\sharp}(\alpha)$
satisfies conditions \eqref{pi-cond1} (which is just the dual of the moment map condition \eqref{mom-2})
and \eqref{pi-cond2}, so Prop.~\ref{prop-inverse} constructs $\pi$ back.

Since the proofs of Propositions \ref{prop-direct} and \ref{prop-inverse} involve long
computations, we will postpone them to the next subsection; we will discuss some examples and implications
of the results first.

\begin{example}(\textit{Nondegenerate quasi-Poisson and quasi-hamiltonian $\mathfrak{g}$-manifolds})\label{ex:nondeg}

It is clear from the correspondence constructed in  Prop.~\ref{prop-direct} that
the singular foliation associated with $\pi$, tangent to $\mathrm{Im}(\pi^{\sharp})+ \mathrm{Im}(\rho_M)\subseteq TM$,
coincides with the singular foliation of the Dirac structure $L$, tangent to $\pr_1(L)\subseteq TM$.
In other words, the Lie algebroids associated with $\pi$ and $L$ have the same leaves, so one is transitive
if and only if the other one is. Note that, for Dirac structures, $\pr_1(L)=TM$ means exactly that
$L$ is defined by a 2-form. As a result, it follows that the correspondence established by
Theorem \ref{quasi-Poisson-hamiltonian} restricts to a one-to-one correspondence between
\textit{nondegenerate} hamiltonian quasi-Poisson manifolds and \textit{presymplectic} realizations 
of Cartan-Dirac structures on Lie groups.
\end{example}

Combining Example \ref{ex:nondeg} with Theorem \ref{quasi-hamiltonian}, we obtain

\begin{corollary} There is a one-to-one correspondence between
\begin{enumerate}
\item[(i)] non-degenerate quasi-Poisson hamiltonian $\mathfrak{g}$-manifolds.
\item[(ii)] quasi-hamiltonian $\mathfrak{g}$-manifolds.
\item[(iii)] presymplectic realizations of $G$ endowed with the Cartan-Dirac structure.
\end{enumerate}
\end{corollary}
Of course, in general, the leaves of a hamiltonian
quasi-Poisson $\mathfrak{g}$-manifold are quasi-hamiltonian $\mathfrak{g}$-manifolds,
which can now be seen as a particular case of Dirac structures having presymplectic
foliations (see also Corollary \ref{cor:restric}). The equivalence of (i) and (ii)
can be found in \cite{AKM}.

The next example answers a question posed in \cite[Ex.~5.2]{SeWe01}.

\begin{example}(\textit{Cartan-Dirac and quasi-Poisson structures on Lie groups})

If $G$ is a Lie group with Lie algebra $\mathfrak{g}$ equipped with a bi-invariant
nondegenerate quadratic form  $\Bi{\cdot,\cdot}$, then one can regard it as a $-\phiG$-twisted Dirac
manifold with respect to the Cartan-Dirac structure $L_G$, or as a hamiltonian
quasi-Poisson $\mathfrak{g}$-manifold. In the latter case,
we consider $G$ acting on itself by conjugation, the quasi-Poisson tensor
is $\pi_G$, defined in Example \ref{ex:quasip}, and the moment map $J:G \to G$
is the identity map, see \cite[Sec.~3]{AKM}.

We claim that these two structures are ``dual'' to each other in the sense of
Theorem \ref{quasi-Poisson-hamiltonian}. Indeed, starting with $\pi_G$ and using
Prop.~\ref{prop-direct}, we see that the corresponding Dirac structure is
$$
L=\{(\pi_G^{\sharp}(\alpha) + \rho(v), C^*(\alpha)+ \sigma(v)),\,|\, (\alpha,v) \in
T^*M \oplus \mathfrak{g}\}.
$$
Since $L_G$ is the image of $(\rho,\sigma)$, it is clear that $L_G \subseteq L$, which
implies that $L_G=L$ since they have the same dimension.
\end{example}

To complete the ``duality'' picture between Dirac realizations of
Cartan-Dirac structures and hamiltonian quasi-Poisson
$\mathfrak{g}$-manifolds, we note that the
correspondence established in Theorem \ref{quasi-Poisson-hamiltonian} preserves maps.

If $(M_i,\pi_i)$ is a hamiltonian quasi-Poisson $\mathfrak{g}$-manifold with moment map
$J_i:M_i\to G$, $i=1,2$, then a map $f:M_1\to M_2$ is a \textbf{hamiltonian quasi-Poisson map}
if $f_*\pi_1=\pi_2$, $f$ is $\mathfrak{g}$-equivariant, and $J_2\circ f = J_1$.
Suppose that $L_i$ is the Dirac structure on $M_i$ corresponding to $\pi_i$, $i=1,2$, via
Theorem \ref{quasi-Poisson-hamiltonian}.

\begin{proposition}\label{prop:maps}
A map $f:(M_1,\pi_1)\to (M_2,\pi_2)$ is a hamiltonian quasi-Poisson map if and only if $f:(M_1,L_1)\to (M_2,L_2)$
is f-Dirac and commutes with the realization maps, $J_2\circ f = J_1$.
\end{proposition}

\begin{proof}
Suppose $f:M_1\to M_2$ is a hamiltonian quasi-Poisson map. In order to check that $f$ is f-Dirac,
we have to compare, at each point of $M_2$, $L_2$ with
\begin{equation}\label{eq:fL1}
f_*L_1=\{(df(X),\beta)\;|\;(X,df^*(\beta))\in L_1\}.
\end{equation}
To simplify the notation, we will denote the infinitesimal actions
of $\mathfrak{g}$ on $M_i$ by $\rho_i$, and
$C_i=1-(1/4)\rho_i\rho^{\vee}dJ$, $i=1,2$. Since $L_2$ corresponds
to $\pi_2$, we have
$$
L_2=\{(\pi_2^{\sharp}(\beta) + \rho_{2}(v), C_2^*(\beta) +
J_2^*\sigma(v))\,|\, \beta \in T^*M_2, v\in \mathfrak{g} \}.
$$
Using that $df\pi_1^{\sharp}df^*=\pi_2^{\sharp}$ (which is another
way of writing  $f_*\pi_1=\pi_2$) and $df\rho_1=\rho_2$ (which is
$f$ $\mathfrak{g}$-equivariance), we deduce that
\begin{equation}\label{eq:map1}
\pi_2^{\sharp}(\beta) +
\rho_2(v)=df(\pi_1^{\sharp}(df^*(\beta))+\rho_1(v)).
\end{equation}
On the other hand, using the $\mathfrak{g}$-equivariance of $f$
and $J_2\circ f=J_1$, it follows that $dfC_1=C_2df$, and we obtain
\begin{equation}\label{eq:map2}
df^*(C_2^*(\beta) +
J_2^*\sigma(v))=C_1^*(df^*(\beta))+J_1^*\sigma(v).
\end{equation}
Since
$$
(\pi_1^{\sharp}(df^*(\beta))+\rho_1(v),
C_1^*(df^*(\beta))+J_1^*\sigma(v)) \in L_1,
$$
for $L_1$ corresponds to $\pi_1$, it follows that, at each point,
$L_2 \subseteq f_*L_1$. But since they have equal dimension, we
conclude that $L_2=f_*L_1$, so $f$ is forward Dirac.

For the converse, suppose that $f_*L_1=L_2$ and $J_2\circ f=J_1$.
It is easy to check that, in this case, $f$ is automatically
$\mathfrak{g}$-equivariant with respect to
\eqref{construct-action}. From that, it follows that
$dfC_1=C_2df$. In order to prove that $f$ is a hamiltonian
quasi-Poisson map, we must still check that $f_*\pi_1=\pi_2$, or,
equivalently, that $df\pi_1^{\sharp}df^*=\pi_2^{\sharp}$. By
Prop.~\ref{prop-inverse}, it suffices to prove that, for $\beta
\in T^*M_2$, $Y=df\pi_1^{\sharp}df^*(\beta) \in TM_2$ satisfies
\begin{equation}\label{eq:cond}
(Y,C_2^*(\beta)) \in L_2\;\; \mbox{ and }\;\;
dJ_2(Y)=-(\rho_2\sigma^{\vee})^*\beta.
\end{equation}
Since $f$ is an f-Dirac map, the first condition in
\eqref{eq:cond} holds since
$$
(\pi_1^{\sharp}df^*(\beta),
df^*(C_2^*(\beta)))=(\pi_1^{\sharp}df^*(\beta),
C_1^*(df^*(\beta))) \in L_1,
$$
for $L_1$ corresponds to $\pi_1$. The second condition holds since
$$
dJ_2(df\pi_1^{\sharp}df^*(\beta))=dJ_1(\pi_1^{\sharp}df^*(\beta))=-(\rho_1\sigma^{\vee})^*df^*(\beta)=
-(\rho_2\sigma^{\vee})^*(\beta).
$$

\end{proof}
We now proceed to the proofs of Propositions \ref{prop-direct} and
\ref{prop-inverse}.

\subsubsection{The proofs}

{\bf Proof of Proposition \ref{prop-direct}:}\\
To simplify our formulas, we set
$$
T= \rho_{M}\rho^{\vee}dJ,
$$ 
and we denote by $\SP{ \cdot, \cdot}$ the pairing between
vector spaces and their duals. 

First, we have to show that $\SP{s(\xi), r(\xi')}$ is antisymmetric in
$\xi, \xi'\in A= T^*M\oplus \mathfrak{g}$. We check this on elements of type $(\alpha, 0)$, $(\beta, 0)$,
and we leave the other cases to the reader. Since $\pi$ is antisymmetric, we only have to show that
$\SP{ T^*\beta, \pi^{\sharp}(\alpha)}$ is antisymmetric in $\alpha$ and $\beta$. 
Using that $dJ\pi^{\sharp}= - (\sigma^{\vee})^{*}(\rho_{M})^{*}$, which is the adjoint
 of the moment map condition (\ref{mom-2}), we see that
\begin{equation}
\label{form0}
\pi(\alpha, T^*\beta)= \SP{T^*\beta, \pi^{\sharp}(\alpha)}= 
-\SP{\rho_{M}^{*}(\beta), \rho^{\vee}(\sigma^{\vee})^{*}\rho_{M}^{*}(\alpha)},
\end{equation}
which is antisymmetric by \eqref{eq:id5}.

We now turn to proving that $s$ satifies (\ref{eq:Dirac}). 
For sections of $A$ of type $\xi= (0, u)$, $\xi'= (0, v)$, with $u, v\in \mathfrak{g}$,
\eqref{eq:Dirac} follows from the next lemma.

\begin{lemma}
\label{lemma-case1}
Given a $\mathfrak{g}$-manifold $M$ and an equivariant map $J: M\rmap G$, then, for any
$u, v\in \mathfrak{g}$,
\[ J^*\sigma([u, v])  = \mathcal{L}_{\rho_{M}(u)}(J^*\sigma(v))-
                                      \mathcal{L}_{\rho_{M}(v)}(J^*\sigma(u))+
                                      d\SP{J^*\sigma(v),\rho_{M}(u)} - 
                                      i_{\rho_{M}(u)\wedge
                                      \rho_{M}(v)}(J^*\phi^{G}).
\]
\end{lemma}

\begin{proof}
Using the equivariance of $J$, $dJ\rho_{M}= \rho$, 
we immediately see that this equation is the pull-back by $J$ of \eqref{eq:Dirac} for 
$\sigma$ (for instance, the last term in the equation equals to $J^*i_{u\wedge v}(\phi^{G})$).
\end{proof}

Next we consider the case $\xi= (0, u)$ and $\xi'= (\alpha, 0)$, which is handled by the next result.
\begin{lemma}
\label{lemma-case2}
Given a $\mathfrak{g}$-manifold $M$, an equivariant map $J: M\rmap G$, and a bivector $\pi$ on $M$ satisfying the moment map
condition (\ref{mom-2}), then, for all $v\in \mathfrak{g}$ and $\alpha\in T^*M$,
\[ 
C^*\Lie_{\rho_M(v)}(\alpha)= 
\mathcal{L}_{\rho_M(v)}(C^*(\alpha))- \mathcal{L}_{\pi^{\sharp}(\alpha)}(J^*\sigma(v))+ 
d\SP{J^*\sigma(v), \pi^{\sharp}(\alpha)}- i_{\rho_M(v)\wedge \pi^{\sharp}(\alpha)}(J^*\phi^{G}). 
\]
\end{lemma}
Although  Lemma \ref{lemma-case2} is more difficult than Lemma \ref{lemma-case1},
it is still simpler than the next case, treated in Lemma \ref{lemma-case3} below.
Since a formula that holds under the same assumptions and is proven by the same method is proven in
detail in Claim \ref{g-invariant} below, we will omit its proof.

Let us now consider $\xi= (\alpha,0 )$ and $\xi'= (\beta, 0)$. Comparing with Lemmas \ref{lemma-case1} and
\ref{lemma-case2}, the greater technical difficulty of this case comes
from the fact that now the formulas involve both $\tri_G$ and $\phi^G$, 
the bracket $[\cdot, \cdot]$ induced by $\pi$ on 1-forms, and require more than just the moment map condition.

\begin{lemma}
\label{lemma-case3}
Given a $\mathfrak{g}$-manifold $M$, an equivariant map $J: M\rmap G$, and an invariant bivector $\pi$ on $M$ satisfying the moment map
condition (\ref{mom-2}), then, for all $\alpha, \beta\in T^*M$,
\[ 
C^*([\alpha, \beta])+ \frac{1}{2}J^*\sigma(i_{(\rho_{M})^*(\alpha\wedge \beta)}\tri_{G}) 
= \mathcal{L}_{\pi^{\sharp}(\alpha)}C^*(\beta)- \mathcal{L}_{\pi^{\sharp}(\beta)}C^*(\alpha)- 
d\SP{C^*(\beta), \pi^{\sharp}(\alpha)}-
i_{\pi^{\sharp}(\alpha)\wedge \pi^{\sharp}(\beta)}(J^*\phi^{G})).
\]
\end{lemma}

\begin{proof}
Since $C= 1- \frac{1}{4}T$, using the definition of $s$ and of the bracket in
$\Gamma(A)$, we see that we can rewrite the equation in the lemma as
\begin{eqnarray}
T^{*}([\alpha, \beta])- \Lie_{\pi^{\sharp}(\alpha)}(T^*(\beta))+ 
\Lie_{\pi^{\sharp}(\beta)}(T^*(\alpha))+ d\pi(\alpha, T^*(\beta)) = &
 2 J^*\sigma i_{\rho_{M}^{*}(\alpha)\rho_{M}^{*}(\beta)}\tri_{G}\nonumber\\ 
&  +4i_{\pi^{\sharp}(\alpha)\wedge\pi^{\sharp}(\beta)}J^*(\phi^{G}). \label{LHSRHS}
\end{eqnarray}
Let us evaluate all the terms of \eqref{LHSRHS} on an arbitrary vector field $X \in \mathcal{X}(M)$.
To simplify the formulas, we set
\begin{equation}\label{eq:notation}
a= \rho_{M}^{*}(\alpha),\;\;  b= \rho_{M}^{*}(\beta), \;\; V= J(X) 
\end{equation}
and we consider the $\mathrm{Hom}(\mathfrak{g}^*, \mathfrak{g})$-valued function on $G$ given by
\[ 
D= \rho^{\vee}(\sigma^{\vee})^{*}.
\]

\begin{claim} \label{claim1} The following formula holds:
\begin{equation}\label{eqclaim1} 
\SP{d\pi(\alpha, T^*(\beta)), X}= \SP{a, \Lie_{V}(Db)}-  \SP{b, \Lie_{V}(Da)}- \SP{a, \Lie_{V}(D)b} .
\end{equation}
\end{claim}

\begin{proof} The left hand side of \eqref{eqclaim1}
is $\Lie_{X}\pi(\alpha, T^*\beta)$. Hence, using (\ref{form0}), it equals
\begin{equation}\label{eqclaim1-2}
- \SP{\Lie_{J(X)}\rho_{M}^{*}(\beta), \rho^{\vee}(\sigma^{\vee})^{*}\rho_{M}^{*}(\alpha)}- 
\SP{\rho_{M}^{*}(\beta), \Lie_{dJ(X)}(\rho^{\vee}(\sigma^{\vee})^{*})}\rho_{M}^{*}(\alpha) .
\end{equation}
With the notation of \eqref{eq:notation}, 
and using $D^*= -D$ (i.e. (\ref{eq:id5})) to rewrite the first term, we see that
\eqref{eqclaim1-2} equals
\[ 
\SP{D\Lie_{V}(b), a}- \SP{b, \Lie_{V}(Da)} .
\]
To obtain \eqref{eqclaim1}, we  write
$D \Lie_{V}(b)= \Lie_{V}(Db)- \Lie_{V}(D)(b)$.
\end{proof}

From the definition of $[\alpha, \beta]$  (\ref{eq:bracket1}), we have
\begin{equation}\label{eqsplit} 
T^*([\alpha, \beta])= T^*\Lie_{\pi^{\sharp}(\alpha)}(\beta)- T^*\Lie_{\pi^{\sharp}(\beta)}(\alpha)
- T^*d\pi(\alpha, \beta)  .
\end{equation}

\begin{claim}\label{claim2}The following formula holds:
\begin{equation}\label{eqclaim2} 
\SP{T^*d\pi(\alpha, \beta), X}= \pi(\Lie_{T(X)}(\alpha), \beta)+ \pi(\alpha, \Lie_{T(X)}(\beta)).
\end{equation}
\end{claim}

\begin{proof}
This follows from the invariance of $\pi$ and the fact that the image of $T$ sits inside that of $\rho_M$.
\end{proof}

Using \eqref{eqsplit} and \eqref{eqclaim2},
 we can split the left hand side of (\ref{LHSRHS}) as 
a difference of two terms which are symmetric to each other. The next claim deals with such a term.

\begin{claim}\label{claim3} The following formula holds:
\begin{eqnarray}
\SP{T^*\Lie_{\pi^{\sharp}(\alpha)}(\beta)- \Lie_{\pi^{\sharp}(\alpha)}(T^*(\beta)), X} &= &
 \pi(\Lie_{T(X)}(\alpha), \beta) + \SP{b, \Lie_{V}(Da)}\nonumber \\
  & & + \SP{b, d\rho^{\vee}( (\sigma^{\vee})^*(a), V)} \label{eqclaim3}
\end{eqnarray}
\end{claim}
\begin{proof}
The left hand side of \eqref{eqclaim3} equals 
\begin{equation}
\label{LHS-claim3}
- \SP{\beta, [\pi^{\sharp}(\alpha), T(X)]+ T([\pi^{\sharp}(\alpha), X])}.
\end{equation}
To rewrite $[\pi^{\sharp}(\alpha), T(X)]$, we note that
\begin{equation}\label{eqclaim3-2}
[\pi^{\sharp}(\alpha), \rho_{M}(\tilde{v})]= -\pi^{\sharp}\Lie_{\rho_{M}(\tilde{v})}(\alpha)+ 
\rho_M\Lie_{\pi^{\sharp}(\alpha)}(\tilde{v}),
\end{equation}
for all $\tilde{v}\in C^{\infty}(M, \mathfrak{g})$: Indeed, 
due to $C^{\infty}(M)$-linearity with respect to $\tilde{v}$, it suffices to check \eqref{eqclaim3-2} for
$\tilde{v}$ constant; in this case,  the equation is just the invariance of $\pi$. 
We now use \eqref{eqclaim3-2} for $\tilde{v}= \rho^{\vee}J(X)$ to get
\[ 
[\pi^{\sharp}(\alpha), T(X)]= -\pi^{\sharp}\Lie_{T(X)}(\alpha)+ 
\rho_M\Lie_{\pi^{\sharp}(\alpha)}(\rho^{\vee}dJ(X)).
\]
We deduce that (\ref{LHS-claim3}) equals to
\begin{align}
&\SP{\beta, \pi^{\sharp}\Lie_{T(X)}(\alpha)}  -  \SP{\rho_{M}^{*}(\beta), 
\Lie_{\pi^{\sharp}(\alpha)}(\rho^{\vee}dJ(X))}+
\SP{\rho_{M}^{*}(\beta), \rho^{\vee}J[\pi^{\sharp}(\alpha), X]}= \nonumber\\
& \pi(\Lie_{T(X)}(\alpha), \beta)- \SP{\rho_{M}^{*}(\beta), \Lie_{\pi^{\sharp}(\alpha)}(\rho^{\vee}J)(X)}.
\label{eqclaim3-3}
\end{align}
On the other hand, for all vector fields $Y$ on $M$ and $\mathfrak{g}$-valued 1-forms $\nu$ on $G$, we have
\[ 
i_X\Lie_{Y}(J^*\nu)= \Lie_{dJ(X)}(\nu(dJ(Y)))+ (d\nu)(dJ(Y), dJ(X)).
\]
Using this identity for $Y= \pi^{\sharp}\alpha$ and $\nu= \rho^{\vee}$ in \eqref{eqclaim3-3}, together
with the dual of the moment map condition, 
$dJ\pi^{\sharp}= - ((\sigma)^{\vee})^{*}(\rho_{M})^{*}$, we obtain
\eqref{eqclaim3} and prove the claim.
\end{proof}

Combining the formulas of claims \ref{claim1}, \ref{claim2} and \ref{claim3}, 
we  conclude that the left hand side of (\ref{LHSRHS}) evaluated at a vector field $X$ is
\begin{equation}
\label{formula-LHS}
\SP{b, d\rho^{\vee} ((\sigma^{\vee})^{*}a, V)}- \SP{a, d\rho^{\vee} ((\sigma^{\vee})^{*}b, V)}- 
\SP{a,\Lie_{V}(D)b}.
\end{equation}
On the other hand, the right hand side of (\ref{LHSRHS}) applied to $X$ equals to
\begin{align}
&2( \tri_{G}(\rho_{M}^{*}(\alpha), \rho_{M}^{*}(\beta), \sigma^*dJ(X))+ 
4\phi^{G}(dJ\pi^{\sharp}(\alpha), dJ\pi^{\sharp}(\beta), dJ(X))) = \nonumber \\  
& 2(\tri_{G}(a, b, \sigma^*V)+ 4\phi^{G}((\sigma^{\vee})^* a, (\sigma^{\vee})^* b, V)), \label{formula-RHS}
\end{align}
where we have used again that $dJ\pi^{\sharp}= - ((\sigma)^{\vee})^{*}(\rho_{M})^{*}$.

To conclude the proof of the lemma,
it suffices to show that $\tri_{G}$ and $\phi^{G}$ are related as follows.

\begin{claim}\label{claim4} 
For all $a, b\in \mathfrak{g}^{*}$ and all vector fields $V$ on $G$, one has
\begin{eqnarray}
\frac{1}{2}\tri_{G}(a, b, \sigma^*V)+ \phi^{G}((\sigma^{\vee})^* a, (\sigma^{\vee})^* b, V) & = &
\frac{1}{4}(-\SP{a, d\rho^{\vee} ((\sigma^{\vee})^{*}b, V)}+\SP{b, d\rho^{\vee} ((\sigma^{\vee})^{*}a, V)}
 \nonumber\\
& &  - \SP{a, \Lie_{V}(D)b}).\label{nustiu}
\end{eqnarray}
\end{claim}

\begin{proof}
It suffices to prove \eqref{nustiu} on elements of type
\[ 
a= u^{\vee}, b= v^{\vee}, V= w_{r},
\]
where  $u,v,w\in \mathfrak{g}$, and
 we recall that $u^{\vee} \in \mathfrak{g}^*$ denotes the dual of $u$ with respect to the quadratic form, and 
$w_{r}$ is the vector field on $G$ obtained from $w$ by right translations. 
We will also denote by $\Ad(u)\in C^{\infty}(G,\mathfrak{g})$
the function $g\mapsto \Ad_{g}(u)$, and we define $\Ad^{-1}(u)$ similarly. 
We will need the explicit formulas for $\sigma^*$ and $(\sigma^{\vee})^*$:
\[ 
\sigma^{*}(w_{r})= \frac{1}{2}(w+ \Ad^{-1}(w))^{\vee},\qquad  
(\sigma^{\vee})^*(u^{\vee})= \frac{1}{2}(u_{r}+ u_{l}).
\]
Using these formulas, combined with 
the invariance of $\tri_{G}$ and $\phi^{G}$, the formula $u_{l}= \Ad(u_{r})$, and the explicit
formulas for $\tri_{G}$ and $\phi^{G}$ on elements of $\mathfrak{g}$, 
one can check that the left hand side of \eqref{nustiu} is
\begin{equation}
\label{firstpart}
\frac{1}{8}(\Bi{[u,v], w+\Ad^{-1}(w)}+ \Bi{[u+\Ad(u), v+\Ad(v)], w}) .
\end{equation}
Since
$\rho^{\vee}$ is the difference between the right and left Maurer-Cartan forms on $G$, we have
\[ 
(d\rho^{\vee})(u_r, v_r)= - [u, v] - \Ad^{-1}[u, v]. 
\]
By the invariance of $\Bi{\cdot,\cdot}$ with respect to $\Ad$, we get that
\begin{eqnarray}
 \SP{a, d\rho^{\vee} ((\sigma^{\vee})^{*}b, V)} & = & 
-\frac{1}{2}(\Bi{u, [v, w]+ \Ad^{-1}([v, w])}+ \Bi{u, [\Ad(v), w]+ [v, \Ad^{-1}(w)]})\nonumber \\
& = &- \frac{1}{2}\Bi{[u+ \Ad(u), v+ \Ad(v)], w}.
\end{eqnarray}
Since $D= \frac{1}{2}(\Ad- \Ad^{-1})$, and $\Lie_{w_r}(D)(v)= -[w, \Ad(v)]- \Ad^{-1}([w, v])$, it follows that
\[ 
\SP{a,\Lie_{V}(D)b}= \frac{1}{2}(\Bi{[u, \Ad(v)], w}+ \Bi{[\Ad(u), v], w}).
\]
Hence the right hand side of \eqref{nustiu} equals to
\[ 
\frac{1}{8}(2\Bi{[u+ \Ad(u), v+ \Ad(v)], w}- \Bi{[u, \Ad(v)], w}- \Bi{[\Ad(u), v], w}),
\]
which is easily seen to coincide with (\ref{firstpart}). This concludes the proof of the claim.
\end{proof}

Using Claim \ref{claim4}, we conclude that \eqref{formula-LHS} and \eqref{formula-RHS} coincide, and this
proves Lemma \ref{lemma-case3}
\end{proof}

From Lemmas \ref{lemma-case1}, \ref{lemma-case2} and \ref{lemma-case3}, it follows that
$s$ is a $-J^*\phi^{G}$-IM form for $A$. To conclude that
$L= \mathrm{Im}(r, s)$ is a Dirac structure, we must still prove that
$L$ has rank $n= \dim(M)$. 
This  follows from the next lemma (which also serves as inspiration
for the proof of Prop.~\ref{prop-inverse}).

\begin{lemma}
\label{lemma-dimension}
The sequence
\[ 0\rmap T^*G \stackrel{j}{\rmap} A\stackrel{(r, s)}{\rmap} L\rmap 0 \]
is exact, where
$j(a)= (-J^*a, \sigma^{\vee}(a))$, $a \in T^*G.$
\end{lemma}

\begin{proof}
The fact that $(r, s)\circ j= 0$ is equivalent to (\ref{mom-2}) and the second formula in (\ref{eq:CC*}).
We define the maps
\begin{equation}
\label{def-U}
U: A\rmap T^*G,\qquad  U(\alpha, v)= -\frac{1}{4}(\rho^{\vee})^{*}\rho_{M}^{*}(\alpha) + \sigma(v) ,
\end{equation}
\begin{equation}
\label{def-i}
i: L\rmap A, \qquad i(X, \alpha)= (\alpha, \frac{1}{4}\rho^{\vee}dJ(X)) .
\end{equation}
We claim that
\begin{equation} \label{identities}
U\circ j= \id, \;\;\; (r, s)\circ i= \id,\;\;  \mbox{ and }\;\; j\circ U+ i\circ (r, s)= \id,
\end{equation}
and these identities imply that the sequence is exact.
For the first identity in \eqref{identities}, write
\[ 
U(j(a))= \frac{1}{4}\rho^{\vee}\rho_{M}^{*}J^*a+ \sigma\sigma^{\vee}a,
\]
and then, using that  $dJ\rho_M= \rho$ and (\ref{eq:id2}), we see that $U\circ j= \id$. The second
identity is immediate from the first and the last ones. To prove the last identity, we
evaluate $j(U(\alpha, v))$: The first component gives
\begin{eqnarray}
 -J^*( -\frac{1}{4}(\rho^{\vee})^{*}\rho_{M}^{*}(\alpha) + \sigma(v)) & = & 
-J^*\sigma(v)+ \alpha- (1- \frac{1}{4}(\rho_M\rho^{\vee}J)^*)\alpha\nonumber\\
 & = & - (J^*\sigma(v)+ C^*(\alpha))+ \alpha
\end{eqnarray}
The second component is
\begin{equation}\label{seccomp}
\sigma^{\vee}( -\frac{1}{4}(\rho^{\vee})^{*}\rho_{M}^{*}\alpha + \sigma(v)).
\end{equation}
But $\sigma^{\vee}(\rho^{\vee})^*\rho_{M}^{*}= - \rho^{\vee}(\sigma^{\vee})^*\rho_{M}^{*}= 
\rho^{\vee}dJ\pi^{\sharp}$,
where we have used (\ref{eq:id5}), and the moment map condition (\ref{mom-2}). Expressing $\sigma^{\vee}\sigma$ 
using (\ref{eq:id1}), we see that \eqref{seccomp} is
\[ 
v- \frac{1}{4}\rho^{\vee}dJ(\rho_{M}(v)+\pi^{\sharp}(\alpha)).
\]
Hence
\[  
j(U(\alpha, v))= (\alpha, v)- (s(\alpha, v), \frac{1}{4}\rho^{\vee}dJ r(\alpha, v) ),
\]
i.e. $j\circ U+ i\circ (r, s)= \id$.
\end{proof}

This concludes the proof of Proposition \ref{prop-direct}.

\noindent{\bf Proof of Proposition \ref{prop-inverse}:}\\
Let $J:(M,L) \to (G,L_G)$ be a Dirac realization.
Identifying $L_G$ with $\mathfrak{g}\ltimes G$, we know that there is an induced action
of $\mathfrak{g}$ on $M$, denoted by $\rho_M$. Spelling out the definition, for $v\in \mathfrak{g}$,
$\rho_{M}(v)$ is the unique vector field satisfying the equations in \eqref{construct-action}.

From the second condition in (\ref{construct-action})
and the fact that $L$ is isotropic, we immediately deduce
\begin{lemma}
\label{lemma-lucky}
For all $(X, \alpha)\in L$,
\begin{equation}
\label{lucky}
\rho_{M}^{*}(\alpha)+ \sigma^*dJ(X) = 0.
\end{equation}
\end{lemma}

Inspired by Lemma \ref{lemma-dimension}, we prove:

\begin{lemma}
\label{lemma-exact}
There is an exact sequence
\[ 0\rmap L\stackrel{i}{\rmap} T^*M\oplus \mathfrak{g}\stackrel{U}{\rmap} T^*G\rmap 0,\]
where $U$ and $i$ are given by (\ref{def-U}) and (\ref{def-i}), respectively.
\end{lemma}

\begin{proof}
First of all, $U$ is surjective since, as in Lemma \ref{lemma-dimension} (and keeping the same notation), 
$U\circ j= \id$.
Next, $U\circ i= 0$ is an immediate consequence of (\ref{eq:id4}) and \eqref{lucky}.
Finally, using the nondegeneracy condition \eqref{eq:nondegreal} for a Dirac realization, 
it follows that $i$ is injective.
By a dimension argument, it follows that the sequence is exact.
\end{proof}

We now concentrate on constructing the quasi-Poisson bivector field $\pi$. Following (ii)
of Prop.~\ref{prop-inverse}, we have

\begin{claim}
$\pi^{\sharp}$ is well defined.
\end{claim}
\begin{proof}
We first show that (\ref{pi-cond1}) and (\ref{pi-cond2}) have a solution $X$, for any given $\alpha$:
the point is that the element
\[ 
(-C^*(\alpha), \frac{1}{4}\rho^{\vee}(\sigma^{\vee})^*\rho_{M}^{*}(\alpha)))
\]
is in the kernel of $U$; this is  a simple computation using (\ref{eq:id2}) and (\ref{eq:id5}).
Hence it must be in the image of $i$. More explicitly, we find that there exists an $X$ such that
\begin{equation}\label{eqinv1} 
(X, C^*(\alpha))\in L,\ \ \rho^{\vee}(dJ(X)+ (\sigma^{\vee})^*\rho_{M}^{*}(\alpha))= 0.
\end{equation}
On the other hand, applying Lemma \ref{lemma-lucky} to
$(X, C^*(\alpha))$, and then using the first equation in (\ref{eq:CC*}) to replace $C\rho_{M}$, we find that
\begin{equation}\label{eqinv2}
\sigma^*(dJ(X)+ (\sigma^{\vee})^*\rho_{M}^{*}(\alpha))= 0.
\end{equation}
Since $\Ker(\sigma^*)\cap \Ker(\rho^{\vee})= 0$, equations \eqref{eqinv1} and \eqref{eqinv2} imply
(\ref{pi-cond1}). The uniqueness of $X$ follows from the nondegeneracy condition \eqref{eq:nondegreal}.
\end{proof}

\begin{claim}
$\pi^{\sharp}$ defines a bivector field $\pi$ which satisfies the moment map condition (\ref{mom-2}).
\end{claim}

\begin{proof}
We have to show that
\[ 
\alpha( \pi^{\sharp}(\beta))+ \beta(\pi^{\sharp}(\alpha))= 0
\]
for all 1-forms $\alpha$ and $\beta$.
Let $X= \pi^{\sharp}(\alpha)$ and $Y= \pi^{\sharp}(\beta)$. Using  (\ref{pi-cond2})
for $(\alpha, X)$ and $(\beta, Y)$, the fact that $L$ is isotropic, and the definition of $C$, we
find that
\begin{equation}\label{eqinv3}
4(\alpha(Y)+ \beta(X))= \alpha(\rho_{M}\rho^{\vee}dJ(Y))+ \beta(\rho_{M}\rho^{\vee}dJ(X)).
\end{equation}
Let us show that the right hand side of \eqref{eqinv3} is zero: using (\ref{pi-cond1}), \eqref{eqinv3} 
becomes:
\[ 
\alpha(\rho_M\rho^{\vee}(\rho_M\sigma^{\vee})^{*}(\beta))+ 
\beta(\rho_M\rho^{\vee}(\rho_M\sigma^{\vee})^{*}(\alpha)),
\]
and this is zero due to (\ref{eq:id5}). On the other hand, (\ref{pi-cond1}) shows that
$dJ\pi^{\sharp}= -(\rho_M\sigma^{\vee})^*$; dualizing it (and using $(\pi^{\sharp})^*= -\pi^{\sharp}$, 
which holds by the first part of the lemma), we obtain the moment map condition.
\end{proof}

\begin{claim}
\label{g-invariant}
The bivector field $\pi$ is $\mathfrak{g}$-invariant.
\end{claim}

\begin{proof} 
We have to show that $\Lie_{\rho_M(v)}(\pi^{\sharp}(\alpha))= \pi^{\sharp}(\Lie_{\rho_M(v)}(\alpha))$ for
$v\in\mathfrak{g}$, and 1-forms $\alpha$. 
For that, it suffices to show that $\Lie_{\rho_M(v)}(\pi^{\sharp}(\alpha))$ satisfies \eqref{pi-cond1} and 
\eqref{pi-cond2}, i.e.,
\begin{align}
& dJ(\Lie_{\rho_M(v)}(\pi^{\sharp}(\alpha)))= -(\rho_M\sigma^{\vee})^*\Lie_{\rho_M(v)}(\alpha), \label{is-equiv-2-a}\\
 & (\Lie_{\rho_M(v)}(\pi^{\sharp}(\alpha)), C^*\Lie_{\rho_M(v)}(\alpha))\in L \label{is-equiv-2}
\end{align}
These conditions are related to Lemma \ref{lemma-case2}. 
Let us first prove (\ref{is-equiv-2}).
Using (\ref{pi-cond2}), (\ref{construct-action}), and the fact that $L$ is isotropic,  we conclude that
\[ 
([\rho_M(v), \pi^{\sharp}(\alpha)], \mathcal{L}_{\rho_M(v)}(C^*\alpha)- 
\mathcal{L}_{\pi^{\sharp}(\alpha)}(J^*\sigma(v))+ d\SP{J^*\sigma(v), \pi^{\sharp}(\alpha)}-
i_{\rho_M(v)\wedge \pi^{\sharp}(\alpha)}(J^*\phi^G)) \in L.
\]
Using Lemma \ref{lemma-case2},
we see that this expression is precisely 
$(\Lie_{\rho_M(v)}(\pi^{\sharp}(\alpha)), C^*\Lie_{\rho_M(v)}(\alpha))$.

Formula (\ref{is-equiv-2-a}) is closely related to the one in Lemma \eqref{lemma-case2}:
the proofs are similar and hold under the same hypothesis (which  might be a bit surprising since
(\ref{is-equiv-2-a}) says that, although the invariance condition on $\pi$ is not assumed, it must be satisfied modulo
the kernel of $J$). Since we have omitted the proof of Lemma \ref{lemma-case2}, we will give the details
for (\ref{is-equiv-2-a}). We evaluate both sides of (\ref{is-equiv-2-a}) on an arbitrary 1-form
$\mu \in \Omega^1(G)$. The left hand side gives 
\begin{eqnarray}
\SP{J^*\mu, [\rho_M(v), \pi^{\sharp}(\alpha)]}& = & 
d(J^*\mu)(\rho_{M}(v), \pi^{\sharp}(\alpha))+ \Lie_{\rho_M(v)}\SP{J^*\mu, \pi^{\sharp}(\alpha)}-
    \Lie_{\pi^{\sharp}(\alpha)}\SP{J^*\mu, \rho_M(v)} \nonumber \\
 & = & - (d\mu)(\rho(v), (\sigma^{\vee})^*\rho_{M}^{*}(\alpha))- 
        \Lie_{\rho_M(v)}\SP{\mu, (\sigma^{\vee})^*\rho_{M}^{*}(\alpha)}+ \nonumber\\
 &&       \Lie_{(\sigma^{\vee})^*\rho_{M}^{*}(\alpha)}\SP{\mu, \rho(v)}\label{expression-late}
\end{eqnarray}
Evaluating $\mu$ on the right hand side, we get
\begin{equation}
 - \SP{\Lie_{\rho_M(v)}(\alpha), \rho_M\sigma^{\vee}\mu}
= - \Lie_{\rho_M(v)}\SP{\alpha, \rho_M\sigma^{\vee}\mu}+ \SP{\alpha, 
[\rho_M(v), \rho_M\sigma^{\vee}\mu]}. \nonumber
\end{equation}
Now, using $[\rho_{M}(v), \rho_{M}(\tilde{v})]= \rho_{M}([v, \tilde{v}])+ \rho_M \Lie_{\rho_M(v)}(\tilde{v})$
for $\tilde{v}= \sigma^{\vee}\mu \in C^{\infty}(M, \mathfrak{g})$, we get
\begin{equation}\label{eqpropin}
- \Lie_{\rho_M(v)}\SP{\rho_{M}^{*}(\alpha), \sigma^{\vee}\mu}+
\SP{\rho_{M}^{*}(\alpha), [v, \sigma^{\vee}\mu]}+ 
\SP{\rho_{M}^{*}(\alpha), \Lie_{\rho_M(v)}(\sigma^{\vee}\mu)}.
\end{equation}
We have to show that this coincides with r.h.s. of (\ref{expression-late}).
Comparing the two formulas, we see that the resulting equation
makes sense for $\rho_{M}^{*}\alpha$ replaced by any element in $C^{\infty}(M, \mathfrak{g}^*)$. 
On the other hand, since the equation is $C^{\infty}(M)$-linear with respect to this element, 
we may assume that the element is a constant
$a\in \mathfrak{g}^*$ (and the remaining appearances of $\rho_M$ become $\rho$). 
The identity to be proven, relating the r.h.s of \eqref{expression-late} and \eqref{eqpropin}, becomes
\begin{equation*}
  -(d\mu)(\rho(v), (\sigma^{\vee})^*a)- \Lie_{\rho(v)}\SP{\mu, (\sigma^{\vee})^*a}+ 
\Lie_{(\sigma^{\vee})^*a}\SP{\mu, \rho(v)}=
   \SP{a, [v, \sigma^{\vee}\mu]},
\end{equation*}
or, equivalently,
\begin{equation}\label{eqequiveq} 
- \mu([\rho(v), (\sigma^{\vee})^*a])= \SP{a, [v, \sigma^{\vee}\mu]}.
\end{equation}
We may assume that $\mu$ is the dual (with respect to the quadratic form) of the vector field 
$w_r$ for some $w\in \mathfrak{g}$, and that
$a$ is the dual of an element $u\in \mathfrak{g}$.  Equation \eqref{eqequiveq} becomes (after multiplying by $2$):
\[ 
- (w_r, [v_r- v_l, u_r+ u_l])= (u, [v, w+ \Ad^{-1}(w)]),
\]
and this can be proven to hold from the invariance of the quadratic form  and the 
identities $[v_l, u_l]= - [v, u]_l$ (see \eqref{eq:sour} for the convention), $[v_r, u_l]= [v_l, u_r]$.
\end{proof}

\begin{claim}
The bivector field $\pi$ is a quasi-Poisson tensor.
\end{claim}

\begin{proof} We must show that $\pi^{\sharp}([\alpha, \beta])= [\pi^{\sharp}(\alpha), \pi^{\sharp}(\beta)]+ 
\frac{1}{2}i_{\alpha\wedge\beta}(\rho_M(\tri_{G}))$.
Using the definition of $\pi^{\sharp}$ ((ii) of Prop.~\ref{prop-inverse}) 
evaluated at $[\alpha, \beta]$, we have to show that
\begin{align}
 & J([\pi^{\sharp}(\alpha), \pi^{\sharp}(\beta)]+ \frac{1}{2}i_{\alpha\wedge\beta}(\rho_M(\tri_{G}))= 
- (\rho_{M}\sigma^{\vee})^*[\alpha, \beta], \\
 & ([\pi^{\sharp}(\alpha), \pi^{\sharp}(\beta)]+ 
   \frac{1}{2}i_{\alpha\wedge\beta}(\rho_M(\tri_{G})), C^*([\alpha, \beta]))\in L. \label{is-quasi-2}
\end{align}
Similarly to the discussion in the previous claim, these conditions are
related to Lemma \ref{lemma-case3}. 
The first equation holds under the same assumptions, and it is proven by the same method, so 
it will be left to the reader (similar to the discussion in the previous proof, the equation tells us that,
although the quasi-Poisson condition is not assumed, it must be satisfied modulo
the kernel of $J$). 

We now prove (\ref{is-quasi-2}).
First, we use that $(\pi^{\sharp}(\alpha), C^*(\alpha))\in L$,  $(\pi^{\sharp}(\beta), C^*(\beta))\in L$,
the fact that $L$ is isotropic, and then apply the formula in Lemma \ref{lemma-case3}, to conclude that
\begin{equation}\label{eqlast1} 
([\pi^{\sharp}(\alpha), \pi^{\sharp}(\beta)], 
C^*([\alpha, \beta])- \frac{1}{2}J^*\sigma i_{\rho_{M}^{*}(\alpha\wedge\beta)}(\tri_{G}))\in L.
\end{equation}
On the other hand, applying the action in (\ref{construct-action}) to 
$v= i_{\rho_{M}^{*}(\alpha\wedge\beta)}(\tri_{G})$ and observing
that $\rho_M(v)= i_{\alpha\wedge\beta}(\rho_M(\tri_{G}))$, we find that
\begin{equation}\label{eqlast2}
(i_{\alpha\wedge\beta}(\rho_M(\tri_{G})), J^*\sigma i_{\rho_{M}^{*}(\alpha\wedge\beta)}(\tri_{G}))\in L.
\end{equation}
Since $L$, at each point, is a vector space, \eqref{eqlast1} and \eqref{eqlast2} imply (\ref{is-quasi-2}).
\end{proof}

\section{Moment maps in Dirac geometry: the global
picture}\label{sec:global}

\subsection{Integrating Lie algebroids and infinitesimal actions}

Lie groupoids are the global counterparts of Lie algebroids. In order to fix our notation, we recall
that a Lie groupoid over a manifold $M$ consists of a manifold $\grd$ together with
surjective submersions $\tar,\sour:\grd \to M$, called \textbf{target} and \textbf{source},
a partially defined multiplication $m: \grd^{(2)}\to \grd$, where
$\grd^{(2)}:=\{(g,h) \in \grd \times \grd\,|\, \sour(g)=\tar(h)\}$, a \textbf{unit section}
$\varepsilon:M \to \grd$ and an \textbf{inversion} $\grd \to \grd$, all related by the
appropriate axioms, see e.g. \cite{SilWein99}. To simplify our notation, we will often
identify an element $x \in M$ with its image $\varepsilon(x) \in \grd$.

For a Lie groupoid $\grd$, the associated Lie algebroid $A(\grd)$
consists of the vector bundle
\begin{equation}\label{eq:sour}
\ker(d\sour)|_M \to M,
\end{equation}
with anchor
$\rho=d\tar: \ker(d\sour)|_M\to TM$ and bracket induced from the
Lie bracket on $\mathcal{X}(\grd)$ via the identification of
sections $\Gamma(\ker(d\sour)|_M)$ with right-invariant vector
fields on $\grd$ tangent to the $\sour$-fibers.

An \textbf{integration} of a Lie algebroid $A$ is a Lie groupoid
$\grd$ together with an isomorphism $A\cong A(\grd)$. Unlike Lie
algebras, not  every Lie algebroid admits an integration, see
\cite{CrFe01} for a description of the obstructions. On the other
hand, if a Lie algebroid is integrable, then there exists a
canonical source-simply-connected integration $\grd(A)$, see
\cite{CrFe01}.

If $M$ is a point, then a Lie groupoid over $M$ is a Lie group, and
the associated Lie algebroid is its Lie algebra.

\begin{example}(\textit{Transformation Lie groupoids})

Let $G$ be a Lie group acting from the left on a manifold $M$.
The associated \textbf{transformation Lie groupoid}, denoted by
$G\ltimes M$, is a Lie groupoid over $M$ with underlying manifold
$G \times M$, source map $\sour(g,x)=x$, target map $\tar(g,x)=g\cdot x$,
and multiplication
$$
(g,x)\cdot (g',x') = (g g', x').
$$
In this case, $A(G\ltimes M)= \mathfrak{g}\ltimes M$, the transformation Lie
algebroid associated with the
infinitesimal action of $\mathfrak{g}$ on $M$ corresponding to
the given $G$-action.
(However, even if a $\mathfrak{g}$-action does not come from a
global action of a Lie group, one can always find a Lie groupoid
integrating the transformation Lie algebroid $\mathfrak{g}\ltimes
M$, see \cite{Da,MoMr}.)
\end{example}

Similarly to infinitesimal actions, Lie groupoids act on maps into
their identity sections: if $\grd$ is a Lie groupoid over $M$,
then a (left) \textbf{action of $\grd$ on a map $J:N \to M$} is a
map $m_N: \grd \times_M N \to N$, $(g,y) \mapsto g\cdot y$,
satisfying
\begin{itemize}
\item[1.] $J(g\cdot y) = \tar(g)$,

\item[2.] $(g g') y = g (g' y)$,

\item[3.] $J(y)\cdot y = y$.
\end{itemize}
Here $\grd \times_M N:= \{(g,y) \in \grd \times N\, |\,
\sour(g)=J(y)\}$. For reasons that will be clear in the next two
subsections, the map $J:N \to M$ is often referred to as the
\textbf{moment map} of the action $m_N$ \cite{MiWe}.

\begin{example}(\textit{Actions of transformation Lie groupoids})\label{ex:acttrangr}

Analogously to Example \ref{ex:actionstransf}, an action $m_N$ of
a transformation Lie groupoid $\grd=G \ltimes M$ on a map $J:N\to
M$ is equivalent to an ordinary action $\overline{m_N}$ of the Lie
group $G$ on $N$ for which $J$ is $G$-equivariant. Indeed, $m_N$
and $\overline{m_N}$ are related by
\begin{equation}
m_N((g,J(y)),y)= \overline{m_N}(g,y), \;\; \mbox{ where } g\in G
\mbox{ and } y\in N.
\end{equation}
\end{example}

The link between infinitesimal and global actions is based on the
following notion: An infinitesimal action $\rho_N$ of a Lie
algebroid $A$ is called \textbf{complete} if $\rho_N(\xi) \in
\mathcal{X}(N)$ is a complete vector field whenever $\xi \in
\Gamma(A)$ has compact support. As in the case of Lie algebras, a
complete action of a Lie algebroid $A$ can be integrated to an
action of its canonical source-simply-connected integration
$\grd(A)$, see e.g. \cite{MoMr}.

\subsection{Poisson maps as moment maps for symplectic groupoid actions}\label{subsec:poissmom}

A 2-form $\omega$ on a Lie groupoid $\grd$ is called
\textbf{multiplicative} if the graph of the groupoid
multiplication $m:\grd^{(2)} \to \grd$ is an isotropic submanifold
of $(\grd,\omega) \times (\grd,\omega) \times (\grd,-\omega)$.
Equivalently, the multiplicativity condition for $\omega$ can be
written as
\begin{equation}\label{eq:compatible}
m^*\omega = \pr_1^*\omega + \pr_2^*\omega,
\end{equation}
where $\pr_i:\grd^{(2)}\to \grd$, $i=1,2$, are the canonical projections.
A \textbf{symplectic groupoid} \cite{We87} is a Lie groupoid together with a multiplicative
symplectic form.

Symplectic groupoids are the global counterparts of Poisson
manifolds in the following sense: If $\pi$ is a Poisson structure
on a manifold $P$ inducing an integrable Lie algebroid structure
on $A=T^*P$ (as in Section \ref{subsec:liealgpoiss}), then the
associated source-simply-connected groupoid $\grd(P):= \grd(A)$
carries a natural multiplicative symplectic structure
\cite{CaFe,CrFe02,MaXu}; on the other hand, on any symplectic
groupoid $(\grd,\omega)$ over a manifold  $P$, condition
\eqref{eq:compatible} automatically implies that $P$ has an
induced Poisson structure uniquely determined by the condition
that the target map $\tar:\grd \to P$ (resp. source map
$\sour:\grd \to P$) is a Poisson map (resp. anti-Poisson map)
\cite{CDW87}.

An \textbf{integration} of a Poisson manifold $(P,\pi)$ is a symplectic groupoid $(\grd,\omega)$
over $P$ for which
the induced Poisson structure coincides with $\pi$. Note that the symplectic form $\omega$ defines a
vector bundle map
\begin{equation}\label{eq:iso}
\ker(d\sour)|_P \longrightarrow  T^*P,\;\; \xi \mapsto
i_{\xi}\omega|_{TP}
\end{equation}
inducing  an isomorphism of Lie algebroids $A(\grd)\cong T^*P$ \cite{CDW87}.
This immediately implies that $\dim(\grd)=2 \dim(P)$.

\begin{example}(\textit{Integrating Lie-Poisson structures})\label{ex:liepoissint}

Let us consider $\mathfrak{g}^*$, equipped with its Lie-Poisson
structure. If $G$ is a Lie group with Lie algebra $\mathfrak{g}$,
then the transformation groupoid $\grd= G\ltimes \mathfrak{g}^*$,
with respect to the coadjoint action, integrates
$T^*\mathfrak{g}^*=\mathfrak{g}\ltimes \mathfrak{g}^*$. The
identification $G \times \mathfrak{g}^* \cong T^*G$ by right
translations induces a multiplicative symplectic form $\omega$ on
$\grd$, in such a way that $(\grd,\omega)$ is a symplectic
groupoid integrating $\mathfrak{g}^*$.
\end{example}

\begin{remark}\label{rem:poisslie2}
The construction of the symplectic groupoid in the previous
example can be extended to the context of Poisson-Lie groups, see
Remark \ref{rem:poisslie1}: If $(G,\pi)$ is a simply-connected
Poisson-Lie group and $G^*$ is its dual, then, assuming that the
dressing action is complete, the transformation groupoid $G\ltimes
G^*$ carries a symplectic structure making it into a symplectic
groupoid integrating $G^*$. (This symplectic structure is
basically the one associated with the semi-direct product Poisson
structure on $G\times G^*$ induced from the action of $G$ on
itself by right multiplication.) For a more general construction
when the actions are not complete, see \cite{LuWe}.
\end{remark}

Let us assume that $P$ is an integrable Poisson manifold.
We have seen that any Poisson map $J:Q \to P$ induces a Lie
algebroid action of $T^*P$ on $Q$. Analogous to the case of Lie algebras, when
this action is complete, it can be ``integrated'' to an action of
$\grd(P)$, the canonical source-simply-connected symplectic
groupoid of $P$. We remark that the completeness of the $T^*P$
action in the Lie algebroid sense coincides with the notion of
$J:Q\to P$ being complete as a Poisson map, i.e., if $f\in
C^{\infty}(P)$ has compact support (or if $X_f$ is complete), then
$X_{J^*(f)}$ is complete.

The global action $m_N:\grd(P)\times_P Q \to Q$ arising in this
way is compatible with the Poisson structure on $Q$ in the sense
that $\gra(m_N)$ is a \textit{lagrangian} submanifold of
$(\grd(P),\pi) \times (Q,\piQ)\times (Q,-\piQ)$\footnote{A
submanifold $C$ of a Poisson manifold $(P,\pi)$ is
\textbf{lagrangian} if, at each $x \in P$, the intersection of
$T_xC$ with $\widetilde{\pi}(T_x^*P)$, the tangent space to the
symplectic leaf at $x$, is a lagrangian subspace of
$\widetilde{\pi}(T_x^*P)$.}, where $\pi$ is the Poisson structure
associated with the symplectic form $\omega$ on $\grd(P)$. Since
inclusions of symplectic leaves of Poisson manifolds are Poisson
maps, an equivalent way to express this compatibility is that the
restricted action $m_N:\grd(P)\times_P S \to S$ to each symplectic
leaf $(S,\omega_S) \hookrightarrow (Q,\piQ)$ satisfies
\begin{equation}
m_N^*\omega_S = \pr_{\st{\grd}}^*\omega + \pr_{\st{S}}^*\omega_S,
\end{equation}
where $\pr_{\st{\grd}}:\grd(P)\times_P S \to \grd(P)$ and $\pr_{\st{S}}:\grd(P)\times_P S \to S$
are the natural projections, see \cite{MiWe,Xu91b}.
On the other hand, if $(Q,\piQ)$
is a Poisson manifold and  $m_N$ is an action of a symplectic groupoid $\grd$
on $J:Q\to P$ compatible with $\pi_Q$ in the sense just described, then $J$ is
automatically a Poisson map (this is just a leafwise version of \cite[Thm.~3.8]{MiWe}).

The next example is the global version of Example
\ref{ex:infinitham}.
\begin{example}(\textit{Global hamiltonian actions})

Consider $\mathfrak{g}^*$ with its Lie-Poisson structure, and let
$G$ be the simply-connected Lie group with Lie algebra
$\mathfrak{g}$. As in Example \ref{ex:infinitham}, the starting
point is a Poisson map $J:Q \to \mathfrak{g}^*$. Note that $J$ is
complete as a Poisson map if and only if the associated
infinitesimal $\mathfrak{g}$-action is by complete vector fields.
In this case, the global action of the symplectic groupoid
$T^*G\cong G\ltimes \mathfrak{g}^*$ is equivalent, in the sense of
Example \ref{ex:actionstransf}, to the hamiltonian $G$-action
obtained by integrating the infinitesimal hamiltonian
$\mathfrak{g}$-action on $Q$.
\end{example}

So, in the previous example, the ``moment'' $J:Q \to
\mathfrak{g}^*$ for the symplectic groupoid action of $T^*G^*$ is
just a momentum map for a hamiltonian $G$-action in the ordinary
sense.

\begin{remark}
Analogously to the previous example and following Remarks
\ref{rem:poisslie1} and \ref{rem:poisslie2}, a Poisson map $J:Q
\to G^*$, where $G^*$ is the dual group to a complete
simply-connected Poisson Lie group, can be ``integrated'' to an
action of the symplectic groupoid $G\ltimes G^*$, which is
equivalent to a $G$-action on $Q$ for which $J$ is equivariant
(with respect to the dressing action on $G^*$). The ``moment'' $J$
in this case coincides with Lu's momentum map \cite{Lu} for a
Poisson action of a Poisson-Lie group on a Poisson manifold.
\end{remark}

\subsection{Dirac realizations as moment maps for presymplectic groupoid actions}\label{subsec:presympgrd}

In order to describe the global actions ``integrating'' Dirac
realizations, we should first identify the global objects
integrating Dirac manifolds, generalizing symplectic groupoids.
This was done in \cite{BCWZ}: if  $\phi$ is a closed 3-form on
$M$, then a \textbf{$\phi$-twisted presymplectic groupoid} over
$M$ is a Lie groupoid $\grd$ over $M$ equipped with a
multiplicative 2-form $\omega$ such that
\begin{itemize}
\item[1.] $d\omega = \sour^*\phi - \tar^*\phi$,

\item[2.] $\dim(\grd)=2\dim(M)$,

\item[3.] $\ker(\omega_x)\cap \ker(d_x\sour) \cap \ker(d_x\tar)
=\{0\}$,  for all $x \in M$.
\end{itemize}
(Twisted presymplectic groupoids are called
\textbf{quasi-symplectic groupoids} in \cite{Xu}.) The multiplicativity of $\omega$ and 
condition 1. in this definition guarantee that the map
\begin{equation}\label{eq:sig}
\sigma_{\omega}: A \to T^*M, \;\;\; \xi \mapsto i_{\xi}\omega|_{TM}
\end{equation}
is a $\phi$-IM form for $A$, while 2. and 3. are
the extra-conditions needed in Lemma \ref{multiplicat} to insure
that the image $L$ of $(\rho, \sigma_{\omega})$ is a
$\phi$-twisted Dirac structure. When $(G, \omega)$ is a symplectic
groupoid, such $L$ is precisely the Dirac structure associated
with the induced Poisson structure on $M$. As proven in
\cite{BCWZ}, $L$ is uniquely determined by the condition that
$\tar$ is an f-Dirac map (resp., $\sour$ is an anti-f-Dirac map).
Conversely, the canonical groupoid $\grd(L)$ integrating the Lie
algebroid associated with a $\phi$-twisted Dirac structure (assuming it is integrable)
is naturally a $\phi$-twisted presymplectic groupoid
\cite[Sec.~5]{BCWZ}. This correspondence generalizes the one
between Poisson manifolds and symplectic groupoids
\cite{CaFe,CrFe02,MaXu} (see also \cite{CaXu} for the integration
of twisted Poisson structures).

We now have all the ingredients to generalize the ``integration''
procedure of Poisson maps to symplectic groupoid actions,
explained  in Section \ref{subsec:poissmom}, to the context of
Dirac geometry. Let $\LM$ be a $\phi$-twisted Dirac structure on $M$
associated with an integrable Lie algebroid.
We call a Dirac realization $J:N\to M$
\textbf{complete} if the induced Lie algebroid action of $\LM$ on
$N$ is complete, in which case it integrates to an action
$m_N:\grd(\LM)\times_M N \to N$, where $(\grd(\LM),\omega)$ is the
canonical twisted presymplectic groupoid associated with $\LM$. In
this situation, we will simply say that the action $m_N$
\textbf{integrates} the realization $J$.

\begin{theorem}\label{thm:globalact}
Let $(M,\LM)$ be a $\phi$-twisted Dirac manifold and assume that $\LM$ is integrable.
A complete Dirac realization $J:N\to M$ integrates to an action $m_N:\grd(\LM)\times_M
N \to N$ satisfying
\begin{equation}\label{eq:compat}
m_N^*\LN = \tau_{\pr_{\st{\grd}}^*\omega}(\pr_{\st{N}}^*\LN),
\end{equation}
where $\pr_{\st{\grd}}$ and $\pr_{\st{N}}$ are the projections
from $\grd(\LN)\times_M N$ onto $\grd(\LN)$ and
$N$, respectively, and $\tau_{\pr_{\st{\grd}}^*\omega}$ denotes a gauge transformation.

Conversely, if $m_N$ is an action of $\grd(\LM)$ on $J:N\to M$
satisfying \eqref{eq:compat}, then $J$ is f-Dirac; if $J$ also
satisfies \eqref{eq:nondegreal}, then it is a Dirac realization
whose integration is $m_N$.
\end{theorem}

In order to prove the theorem, we need the following result.

\begin{lemma}\label{lem:globalact}
Let $(M,\LM)$ be a $\phi$-twisted Dirac manifold and assume that $\LM$ is integrable.
Let $m_N:\grd(\LM)\times_M N \to N$ be an action of $\grd(\LM)$ on $J:N\to M$, and
assume that $N$ is equipped with a $J^*\phi$-twisted presymplectic form $\omega_N$.
Then $J$ is an f-Dirac map if and only if
\begin{equation}\label{eq:compat2}
m_N^*\omega_N = \pr_{\st{N}}^*\omega_N + \pr_{\st{\grd}}^*\omega.
\end{equation}
\end{lemma}
\begin{proof}
To simplify the notation, let $\grd = \grd(\LM)$, and let us denote by $A$ the corresponding Lie algebroid
(which is just $\LM$). The source and target maps in $\grd$ are denoted by $\sour$ and $\tar$.
Also, let
$\omega_1 = m_N^*\omega_N - \pr_{\st{N}}^*\omega_N$ and
$\omega_2=\pr_{\st{\grd}}^*\omega$. With these definitions, our goal is to show
that $J$ is f-Dirac if and only if $\omega_1=\omega_2$.

The key observation is that if we regard
$\grd\times_M N$ as a transformation Lie groupoid over $N$, with source $\pr_{\st{N}}$ and target $m_N$,
a direct computation shows that both $\omega_1$ and $\omega_2$ are multiplicative. Hence, by
\cite[Thm.~2.5]{BCWZ}, $\omega_1=\omega_2$ if and only if the corresponding bundle maps
$$
\sigma_{\omega_i}: A\times_M N \to T^*N,\;\; \xi_y \mapsto \sigma_{\omega_i}=(i_{\xi_y}\omega_i)|_{TN}
$$
$i=1,2$, see \eqref{eq:sig}, coincide. For $\xi_y \in A\times_M N$ ($\xi \in A_x$ and $x=J(y)$) and $Y \in TN$
(as usual, we identity $TN$ with $T\varepsilon (N)$, where $\varepsilon:N \to \grd \times_M N$ is the identity
section), we have
\begin{equation}\label{eq:sigma1and2}
\sigma_{\omega_1}(\xi_y,Y)=\omega_N(dm_N(\xi_y),Y)\;\; \mbox{ and } \;\;
\sigma_{\omega_2}(\xi_y,Y)=\omega(\xi,dJ(Y)).
\end{equation}
For the first identity in \eqref{eq:sigma1and2}, we used that $i_{\xi_y}\pr_N^*\omega_N=0$ for
$\xi_y \in A\times_M N$, since $\pr_N$ is the source map in $\grd\times_M N$, and $A\times_MN$
is its Lie algebroid, which is tangent to the source fibres along the identity section.

Since $\LM=\{(d\tar(\xi),i_{\xi}\omega|_{TM})\;|\; \xi \in A\}$, $J:N\to M$ being f-Dirac
means that
\begin{equation}\label{eq:Jfdir}
\{(d\tar(\xi),i_{\xi}\omega|_{TM})\;|\; \xi \in A\} = \{(dJ(Y),\alpha)\;|\; i_Y\omega_N=J^*\alpha\}.
\end{equation}
But, for $\xi_y \in A\times_M N$, we have $dJ(dm_N(\xi_y))=d\tar(\xi)$. It then follows from \eqref{eq:Jfdir} that
$$
\omega(\xi,dJ(Y))=\omega_N(dm_N(\xi_y),Y)
$$
for all $Y \in TN$, which implies that  $\sigma_{\omega_1} = \sigma_{\omega_2}$, i.e., $\omega_1=\omega_2$.

The converse follows from the same arguments, reversing the steps.
\end{proof}

We can now prove Theorem \ref{thm:globalact}:

\begin{proof}
We keep writing $\grd$ for $\grd(\LM)$.
Suppose that $m_N$ integrates a Dirac realization $J:N\to M$.
The bundles $m_N^*\LN$ and $\pr_{\st{N}}^*\LN$, seen as subbundles of
$T(\grd\times_M N)\oplus T^*(\grd\times_M N)$, have the same
projection onto the first factor: at a point $(g,y)$, they both coincide
with $T_g\grd \times T_y\mathcal{O}$, where $\mathcal{O}$ is the leaf of $\LN$
through $y$. Note that, since $\pr_{\st{N}}$ is a submersion, $\pr_{\st{N}}^*\LN$
is a smooth subbundle, so it is an honest Dirac structure.

By Corollary \ref{cor:restric}, since $J$ is a Dirac realization of $M$, its restriction to any leaf
of $\LN$, $(\mathcal{O},\theta)$, is a presymplectic realization,
 and $m_N$ is tangent to the leaves. By Lemma \ref{lem:globalact},
$$
m_N^*\theta = \pr_{\st{N}}^*\theta + \pr_{\st{\grd}}^*\omega,
$$
which implies the compatibility \eqref{eq:compat}.

Conversely, \eqref{eq:compat} implies that $J$ is tangent to the leaves of $\LN$. Restricting
$m_N$ to these leaves, \eqref{eq:compat} amounts to \eqref{eq:compat2}. So, by Lemma \ref{lem:globalact},
$J$ is an f-Dirac map when restricted to each leaf, which implies that $J$ is f-Dirac by Corollary \ref{cor:incl}.
The last statement follows from Corollary \ref{cor:restric} and a direct check.
\end{proof}

\begin{remark}

The presymplectic groupoid actions resulting from
\textit{presymplectic} realizations are exactly the ``modules''
considered in the Morita theory developed in \cite{Xu} to compare
various notions of moment maps.  More general Dirac realizations
give rise to more general ``hamiltonian spaces'' which still fit
with the constructions in \cite{Xu}.
\end{remark}

We will discuss examples of the ``integration'' in Theorem
\ref{thm:globalact} related to ``quasi'' hamiltonian actions in
Section \ref{subsec:AMMquasi}.
\subsection{Reduction in Dirac geometry}
Just as in Poisson geometry, one can also carry out reduction in
the context of Dirac manifolds. The general construction described
in this section recovers reduction procedures in various settings, including
\cite{AKM,AMM,MiWe,Xu}.

The set-up is as follows. Let $J:N \to M$ be a Dirac realization
of a $\phi$-twisted Dirac manifold $(M,\LM)$. Let $x \in M$ be a
regular value of $J$, and consider the submanifold $\iota:
\mathcal{C}=J^{-1}(x) \hookrightarrow N$. Following
\cite{MiWe,Xu}, let $\mathfrak{l}_x = \ker(\rho)_x$ be the
isotropy Lie algebra of $\LM$ at $x$. Since the anchor $\rho$ is
the projection $\pr_1|_{\LM}$, it follows that
\begin{equation}
\mathfrak{l}_x = (\LM \cap T^*M)_x.
\end{equation}
The induced Lie algebroid action of $\LM$ on $J:N\to M$ defines a
vector bundle morphism $\LM\times_M N \to TN$, and a simple
computation shows that this morphism gives rise to an action of
the Lie algebra $\mathfrak{l}_x$ on $\mathcal{C}$. Our object of
interest is the orbit space $\mathcal{C}/\mathfrak{l}_x$.

\begin{lemma}\label{lem:reduction}
If the stabilizer algebras of the $\mathfrak{l}_x$-action on $\mathcal{C}$
have constant dimension (on each component), then
$\iota^*\LN$ is a (untwisted) Dirac structure on $\mathcal{C}$.
\end{lemma}
\begin{proof}
As mentioned in Section \ref{subsec:diracmaps}, the conclusion in the lemma holds as
long as we show that $\iota^*\LN$ is a \textit{smooth} subbundle
of $T\mathcal{C}\oplus T^*\mathcal{C}$.

As a vector bundle, $\iota^*\LN$ is naturally identified with
\begin{equation}\label{eq:frac}
\frac{\LN\cap(T\mathcal{C}\oplus T^*N)}{(\LN\cap
T\mathcal{C}^\circ)},
\end{equation}
see \cite{Cou90}, and $\iota^*\LM$ will be smooth if we show that both bundles in
\eqref{eq:frac} are smooth. For that, it
suffices to show that each one has constant dimension. But since
their quotient $\iota^*\LN$ has constant dimension, it suffices to show that
either $\LN\cap(T\mathcal{C}\oplus T^*N)$ or $\LN\cap
T\mathcal{C}^\circ$ has constant dimension. We will prove that for $\LN\cap T\mathcal{C}^\circ$.

On one hand,
$$
\LN\cap T\mathcal{C}^\circ = \{(0,\beta) \in \LN\;|\;
\iota^*\beta=0\} = \{(0, dJ^*\alpha) \in \LN\;|\; \alpha \in
T^*M\}.
$$
It follows from $J$ being an f-Dirac map that if $(0,dJ^*\alpha)
\in \LN$, then $(0,\alpha) \in \LM$. So, if $\rho_N$ is the
infinitesimal action of $\LM$ on $J$, we can write
$$
\LN\cap T\mathcal{C}^\circ = \{(0, dJ^*\alpha) \in \LN\;|\; \alpha
\in T^*M\} \cong \frac{\ker(\rho_N)\cap \LM \cap
T^*M}{\ker(dJ^*)}.
$$
But $\ker(\rho_N)\cap \LM \cap T^*M$ is the stabilizer of the $\mathfrak{l}_x$-action
on $\mathcal{C}$, which is assumed to have constant dimension. Since $x$ is a regular value,
$dJ$ has maximal rank on $\mathcal{C}$, so $\ker(dJ^*)$ also has constant dimension.
As a result, the dimension of \eqref{eq:frac} is constant, and $\iota^*\LN$ is a
smooth bundle.

Finally, note that $\iota^*\LN$ is a $(\iota^*J^*\phi)$-Dirac
structure on $\mathcal{C}$, but $\iota^*J^*\phi=0$. So
$\iota^*\LN$ is an ordinary Dirac structure.
\end{proof}

We now show that the quotient $\mathcal{C}/\mathfrak{l}_x$ carries
a natural Poisson structure.

\begin{theorem}\label{thm:reduction}
Suppose that the orbit space $\mathcal{C}/\mathfrak{l}_x$ is a smooth manifold so that
projection $\mathcal{C} \rmap \mathcal{C}/\mathfrak{l}_x$ is a submersion.
Then there is a
unique Poisson structure $\pi_{\st{red}}$ on
$\mathcal{C}/\mathfrak{l}_x$ for which the projection
$(\mathcal{C},\iota^*\LN) \to
(\mathcal{C}/\mathfrak{l}_x,\pi_{\st{red}})$ is an f-Dirac map.
\end{theorem}

\begin{remark}\label{rem:bandf}
The projection $(\mathcal{C},\iota^*\LN) \to
(\mathcal{C}/\mathfrak{l}_x,\pi_{\st{red}})$ is also a b-Dirac map, and this
property characterizes $\pi_{\st{red}}$ uniquely as well.
\end{remark}

\begin{proof}
It follows from our assumptions that
the $\mathfrak{l}_x$-orbits on $\mathcal{C}$ have constant dimension, so the same holds
for the stabilizer algebras.
By Lemma \ref{lem:reduction}, $(\mathcal{C},\iota^*\LN)$ is a
Dirac manifold.

The \textit{admissible functions} on
$(\mathcal{C},\iota^*\LN)$, i.e., the set of functions on
$\mathcal{C}$ whose differential vanish on
$\ker(\iota^*\LN)=\iota^*\LN \cap T\mathcal{C}$ form a Poisson
algebra, see \cite[Sec.~2.5]{Cou90}, under the bracket
$$
\{f,g\}:=\Lie_{X_f}g,
$$
where $X_f$ is a local vector field such that $(X_f,df) \in
\iota^*\LN$. We will show that this Poisson algebra induces a
Poisson structure on $\mathcal{C}/\mathfrak{l}_x$ by showing that
the kernel of $\iota^*\LN$ coincides with the
$\mathfrak{l}_x$-orbits, i.e.,
\begin{equation}\label{eq:reduc}
\ker(\iota^*\LN) =\rho_N(\mathfrak{l}_x).
\end{equation}
On one hand,
\begin{eqnarray*}
\iota^*\LN\cap TQ&=& \{Y\in TQ\;|\; \exists \beta \in T^*N \mbox{
with } (Y,\beta) \in \LN, \; \iota^*\beta=0\}\\
&=& \{Y\in TQ\;|\; \exists \alpha \in T^*M \mbox{ with }
(Y,dJ^*\alpha)\in \LN\}.
\end{eqnarray*}
But since $J$ is f-Dirac and $dJ(Y)=0$, we can write
$$
\iota^*\LN\cap TQ= \{Y\in TQ\;|\; \exists \alpha \in \LM\cap T^*M
\mbox{ with } (Y,dJ^*\alpha)\in \LN\}.
$$
On the other hand,
\begin{eqnarray*}
\rho_N(\mathfrak{l}_x)&=&\{Y \in TN \;|\;\exists \alpha \in
\LM\cap T^*M \mbox{ with } (Y,dJ^*\alpha) \in \LN,\; dJ(Y)=0\}\\
&=& \{Y\in TQ\;|\; \exists \alpha \in \LM\cap T^*M \mbox{ with }
(Y,dJ^*\alpha)\in \LN\}.
\end{eqnarray*}
So \eqref{eq:reduc} follows.

The fact that the projection $(\mathcal{C},\iota^*\LN) \to
(\mathcal{C}/\mathfrak{l}_x,\pi_{\st{red}})$ is an f-Dirac map and the
claim in Remark \ref{rem:bandf} follow from 
a direct computation, see e.g. \cite{BuRa02}. 
\end{proof}
Of course, if the Dirac realization $J:N\to M$ is complete, one
can state Theorem \ref{thm:reduction} in terms of the action of
the isotropy group of $\grd(\LM)$ at $x \in M$ on $\mathcal{C} =
J^{-1}(x)$. Versions of Theorem \ref{thm:reduction} can also be derived
when this action is locally free and the quotient is an orbifold, as well
as for more general ``intertwiner spaces'' in the sense of  \cite{Xu}.

\begin{remark}(\textit{Other reductions})\label{rem:redex}

The following are important particular cases of the reduction
in Theorem \ref{thm:reduction}:
\begin{itemize}
\item[-] If $M$ is Poisson and $J:N \to M$
is a symplectic realization,  we  recover \cite[Thm.3.12]{MiWe};
in particular, when $M=\mathfrak{g}^*$, this reduces to
Marsden-Weinstein classical theorem \cite{MW}, and when $M=G^*$,
the dual of a Poisson-Lie group, we get Lu's reduction \cite{Lu}.
If $J:N \to M$ is a Poisson map, we get the ``Poisson-version''
of these results.

\item[-] If
$M$ is $\phi$-twisted Dirac and $J:N\to M$ is a presymplectic
realization, then we obtain Xu's reduction \cite[Thm.~3.17]{Xu};
in particular, when $M$ is a Lie group equipped with Cartan-Dirac structure,
one recovers the quasi-hamiltonian reduction of \cite{AMM}.

\item[-] If $J:N \to G$ is a general Dirac realization of a Lie group with Cartan-Dirac structure,
then we recover the reduction of quasi-Poisson manifolds of \cite{AKM} via the identification
established in Theorem \ref{quasi-Poisson-hamiltonian}, see Remark
\ref{rem:qred} below.
\end{itemize}
\end{remark}

\subsection{AMM-groupoids and hamiltonian quasi-Poisson $G$-manifolds}\label{subsec:AMMquasi}
We now discuss global actions, in the sense of Theorem \ref{thm:globalact},
associated with complete Dirac realizations of Cartan-Dirac structures.

Let $G$ be a Lie groups equipped with a Cartan-Dirac structure $L_G$
with respect to a bi-invariant nondegenerate quadratic form $\Bi{\cdot,\cdot}$. The first step is to identify
$\grd(L_G)$, the canonical presymplectic groupoid integrating $L_G$.

As shown in \cite[Sec.~7]{BCWZ}, $\grd(L_G)$ is closely related to the
AMM-groupoids of \cite{BeXuZh}:
if $\grd = G\ltimes G$ is the transformation groupoid with
respect to the conjugation action, then the 2-form \cite{AMM}
$$
\omega_{(g,x)}=\frac{1}{2}(\Bi{\Ad_xp_g^*\lambda,p_g^*\lambda} +
\Bi{p_g^*\lambda,p_x^*(\lambda+\overline{\lambda})}),
$$
where $p_g,p_x:G\times G \to G$ are the first and second
projections, and $\lambda$ and $\overline{\lambda}$ are the left
and right Maurer-Cartan forms, makes $\grd$ into a $\phi^G$-twisted
presymplectic groupoid. If $G$ is simply-connected, then
$(\grd,\omega)$ is isomorphic to  $\grd(L_G)$, the canonical
source-simply-connected integration of $L_G$. In general,
$\grd(L_G)$ is obtained from the AMM
groupoid by pulling back $\omega$ to $\widetilde{G}\ltimes G$,
where $\widetilde{G}$ is the universal cover of $G$
\cite[Thm.~7.6]{BCWZ}.
As a result, just as Lie-Poisson structures
``integrate'' to cotangent bundles of Lie groups, see Example \ref{ex:liepoissint}, Cartan-Dirac
structures ``integrate'' to the ``double'' $(G\times G,\omega)$ in
the sense of \cite{AMM}.

For simplicity, let $G$ be simply connected. A complete Dirac realization $J:M\to G$
induces a presymplectic groupoid action of $(\grd,\omega)$, as in Theorem \ref{thm:globalact},
which is equivalent to a $G$-action on $M$ for which $J$ is $G$-equivariant, see Example
\ref{ex:acttrangr}; this $G$-action is just an integration of the infinitesimal $\mathfrak{g}$-action
which makes $M$ into a quasi-Poisson $\mathfrak{g}$-manifold, as constructed in Proposition \ref{prop-inverse}.
So $M$ becomes a hamiltonian quasi-Poisson $G$-manifold for which $J:M \to G$
is the group valued moment map \cite{AKM}. This construction yields the following global version
of Theorem \ref{quasi-Poisson-hamiltonian}.

\begin{theorem}\label{quasi-Poisson-hamiltonian-global}
There is a one-to-one correspondence between complete Dirac realizations of $(G,L_G)$
and hamiltonian quasi-Poisson $G$-manifolds.
\end{theorem}

\begin{corollary}
There is a one-to-one correspondence between compact Dirac realizations of $(G,L_G)$
and compact hamiltonian quasi-Poisson $G$-manifolds.
\end{corollary}

Of course, a global version of Prop.~\ref{prop:maps} also holds.

\begin{remark}(\textit{Reduction})\label{rem:qred}

Given a Dirac realization of $(G,L_G)$, $J:M\to G$, the Dirac reduction  of
Theorem \ref{thm:reduction} produces Poisson spaces $J^{-1}(g)/G_g$, where $G_g$
is the centralizer of $g\in G$.
Using Remark \ref{rem:redex} and \cite[Prop.~10.6]{AKM}, one can check
that these are the same Poisson spaces obtained by quasi-Poisson reduction
\cite[Thm.6.1]{AKM} if we regard $M$ as a hamiltonian quasi-Poisson $G$-manifold
instead.
\end{remark}

\begin{footnotesize}

\end{footnotesize}

\end{document}